\input ./style/arxiv-general.cfg
\documentclass[aos,MSNbibl,secthm,seceqn,nameyear,dvips]{arximspdf}
\makeatletter
   \@ifpackageloaded{graphicx}{}{\usepackage{graphicx}}
\makeatother

\doi{10.1214/15-AOS1379}
\volume{44}
\issue{5}
\pubyear{2016}
\firstpage{1957}
\lastpage{1987}
\docsubty{FLA}

\makeatletter
\newcommand{\rrvert}{\vert}
\newcommand{\llvert}{\vert}
\newtheorem{lem}{Lemma}
\newtheorem{cor}{Corollary}[section]
\newproclaim{rem}{Remark}[section]
\newproclaim{exa}{Example}[section]
\newproclaim{defn}{Definition}[section]
\newproclaim{ass}{Assumption}[section]
\makeatother

\begin{document}
\begin{frontmatter}

\title{Estimation in nonlinear regression with Harris recurrent Markov chains\thanksref{T1}}
\runtitle{Nonlinear and nonstationary regression}

\begin{aug}
\author[A]{\fnms{Degui}~\snm{Li}\thanksref{T2}\ead[label=e1]{degui.li@york.ac.uk}},
\author[B]{\fnms{Dag}~\snm{Tj{\o}stheim}\ead[label=e2]{Dag.Tjostheim@math.uib.no}}
\and
\author[C]{\fnms{Jiti}~\snm{Gao}\thanksref{T3}\corref{}\ead[label=e3]{jiti.gao@monash.edu}}
\runauthor{D. Li, D. Tj{\o}stheim and J. Gao}
\affiliation{University of York, University of Bergen and Monash University}
\address[A]{D. Li\\
Department of Mathematics\\
University of York\\
Heslington Campus\\
York, YO10 5DD\\
United Kingdom\\
\printead{e1}}
\address[B]{D. Tj{\o}stheim\\
Department of Mathematics\\
University of Bergen\\
Post box 7800\\
5020 Bergen\\
Norway\\
\printead{e2}}
\address[C]{J. Gao\\
Department of Econometrics and\\
\quad Business Statistics\\
Monash University at Caulfield\\
Caulfield East, Victoria 3145\\
Australia\\
\printead{e3}}
\end{aug}
\thankstext{T1}{Supported by the Norwegian Research Council.}
\thankstext{T2}{Supported in part by an Australian Research Council Discovery Early Career Researcher Award
(DE120101130).}
\thankstext{T3}{Supported by two Australian Research Council Discovery Grants under Grant Numbers: DP130104229
and DP150101012.}

\footnotetext{\textbf{Tribute}: While this paper was in the process of being published, we
heard that Professor Peter Hall, one of the most significant contributors to the areas
of nonlinear regression and time series analysis, sadly passed away. The fundamental
work done by Professor Peter Hall in the area of martingale theory,
represented by the book (with Christopher C. Heyde):
Hall, P. and Heyde, C. [\textit{Martingale Limit Theory and Its Applications} (1980) Academic Press],
enables the authors of this paper in using martingale theory as an
important tool in dealing with all different types of estimation and
testing issues in econometrics and statistics. In a related
{\em Annals} paper by Gao, King, Lu and Tj{\o}stheim
[\textit{Ann. Statist.} \textbf{37} (2009) 3893--3928],
Theorem 3.4 of Hall and Heyde (1980) plays an essential role
in the establishment of an important theorem. In short, we
would like to thank the \mbox{Co-}Editors for including our paper in this
dedicated issue in honour of Professor Peter Hall's fundamental contributions to statistics and theoretical econometrics.}

%
\received{\smonth{1} \syear{2014}}
%
\revised{\smonth{8} \syear{2015}}

%
\begin{abstract}
In this paper, we study parametric nonlinear regression under the
Harris recurrent Markov chain framework. We first consider the
nonlinear least squares estimators of the parameters in the
homoskedastic case, and establish asymptotic theory for the proposed
estimators. Our results show that the convergence rates for the
estimators rely not only on the properties of the nonlinear regression
function, but also on the number of regenerations for the~Harris
recurrent Markov chain. Furthermore, we discuss the estimation of the
parameter vector in a conditional volatility function, and apply our
results to the nonlinear regression with $I(1)$ processes and derive an
asymptotic distribution theory which is comparable to that obtained by
Park and Phillips [\textit{Econometrica} \textbf{69} (2001) 117--161].
Some numerical studies including simulation
and empirical application are provided to examine the finite sample
performance of the proposed approaches and results.
\end{abstract}\vspace*{15pt}

%
\end{frontmatter}

\section{Introduction}\label{sec1}

In this paper, we consider a parametric nonlinear regression model
defined by
%
\begin{eqnarray}
Y_t & = & g(X_t,\theta_{01},
\theta_{02},\ldots,\theta_{0d})+e_t
\nonumber
\\[-8pt]
\label{eq1.1}
\\[-8pt]
\nonumber
& =: & g(X_t,{\bolds\theta}_0)+e_t,
\qquad t=1,2, \ldots, n,
\end{eqnarray}
where ${\bolds\theta}_0$ is the true value of the $d$-dimensional
parameter vector such that
\[
{\bolds\theta}_0=(\theta_{01},\theta_{02},
\ldots,\theta _{0d})^\tau \in\Theta\subset{
\mathbb{R}}^d
\]
and $g(\cdot,\cdot): {\mathbb{R}}^{d+1}\rightarrow{\mathbb{R}}$ is
assumed to be known. Throughout this paper, we assume that $\Theta$ is
a compact set and ${\bolds\theta}_0$ lies in the interior of
$\Theta$, which is a standard assumption in the literature. How to
construct a consistent estimator for the parameter vector ${\bolds\theta}_0$ and derive an asymptotic theory are important issues in
modern statistics and econometrics. When the observations $(Y_t, X_t)$
satisfy stationarity and weak dependence conditions, there is an
extensive literature on the theoretical analysis and empirical
application of the above parametric nonlinear model and its extension;
see, for example, \citet{J69}, \citet{M70} and \citet{W81} for
some early references, and \citet{SW92}, \citet{L94},
\citet{S00} and \citet{LN08} for recent relevant works.

As pointed out in the literature, assuming stationarity is too
restrictive and unrealistic in many practical applications. When
tackling economic and financial issues from a time perspective, we
often deal with nonstationary components. For instance, neither the
consumer price index nor the share price index, nor the exchange rates
constitute a stationary process. A traditional method to handle such
data is to take the first-order difference to eliminate possible
stochastic or deterministic trends involved in the data, and then do
the estimation for a stationary model. However, such differencing may
lead to loss of useful information. Thus, the development of a modeling
technique that takes both nonstationary and nonlinear phenomena into
account in time series analysis is crucial. Without taking differences,
\citet{PP01} (hereafter PP) study the nonlinear regression~(\ref{eq1.1}) with the regressor $\{X_t\}$ satisfying a unit root [or
$I(1)$] structure, and prove that the rates of convergence of the
nonlinear least squares (NLS) estimator of ${\bolds\theta}_0$
depend on the properties of $g(\cdot,\cdot)$. For an integrable
$g(\cdot,\cdot)$, the rate of convergence is as slow as $n^{1/4}$, and for an
asymptotically homogeneous $g(\cdot,\cdot)$, the rate of convergence
can achieve the $\sqrt{n}$-rate and even $n$-rate of convergence. More
recently, \citet{CW12} consider the same model structure as
proposed in the PP paper and then establish some corresponding results
under certain technical conditions which are weaker than those used in
the PP paper.

As also pointed out in a recent paper by \citet{MKT12},
the null recurrent Markov process is a nonlinear generalization of the
linear unit root process, and thus provides a more flexible framework
in data analysis. For example, \citet{GTY13} show that the
exchange rates between British pound and US dollar over the time period
between January 1988 and February 2011 are nonstationary but do not
necessarily follow a linear unit root process [see also \citet{BRS08} for a similar discussion of the exchange rates between
French franc and German mark over the time period between December 1972
and April 1988]. Hence, \citet{GTY13} suggest using the
nonlinear threshold autoregressive (TAR) with stationary and unit root
regimes, which can be proved as a $1/2$-null recurrent Markov process;
see, for example, Example~\ref{exa21} in Section~\ref{sec2.2} and
Example~\ref{exa61}
in the empirical application (Section~\ref{sec6}).

Under the framework of null recurrent Markov chains, there has been an
extensive literature on nonparametric and semiparametric estimation
[\citet{KT01}, \citeauthor{KMT07} (\citeyear{KMT07},
\citeyear{KMT10}),
\citet{LLC09}, \citet{S11},
\citet{CGL12},
\citet{GKLT14}], by using the technique of the split chain
[\citet{N84}, \citet{MT09}], and the generalized ergodic
theorem and functional limit theorem developed in
\citet{KT01}. As far as we know, however, there is virtually no work
on the parametric estimation of the nonlinear regression model (\ref{eq1.1}) when the regressor $\{X_t\}$ is generated by a class of Harris
recurrent Markov processes that includes both stationary and
nonstationary cases. This paper aims to fill this gap. If the function
$g(\cdot,\cdot)$ is integrable, we can directly use some existing
results for functions of Harris recurrent Markov processes to develop
an asymptotic theory for the estimator of ${\bolds\theta}_0$. The
case that $g(\cdot,\cdot)$ belongs to a class of asymptotically
homogeneous functions is much more challenging, as in this case the
function $g(\cdot,\cdot)$ is no longer bounded. In nonparametric or
semiparametric estimation theory, we do not have such problems because
the kernel function is usually assumed to be bounded and has a compact
support. Unfortunately, most of the existing results for the asymptotic
theory of the null recurrent Markov process focus on the case where
$g(\cdot,\cdot)$ is bounded and integrable [c.f., \citeauthor{C99} (\citeyear{C99}, \citeyear{C00})].
Hence, in this paper, we first modify the conventional NLS estimator
for the asymptotically homogeneous $g(\cdot,\cdot)$, and then use a
novel method to establish asymptotic distribution as well as rates of
convergence for the modified parametric estimator. Our results show
that the rates of convergence for the parameter vector in nonlinear
cointegrating models rely not only on the properties of the function
$g(\cdot,\cdot)$, but also on the magnitude of the regeneration number
for the null recurrent Markov chain.

In addition, we also study two important issues, which are closely
related to nonlinear mean regression with Harris recurrent Markov
chains. The first one is to study the estimation of the parameter
vector in a conditional volatility function and its asymptotic theory.
As the estimation method is based on the log-transformation, the rates
of convergence for the proposed estimator would depend on the property
of the log-transformed volatility function and its derivatives.
Meanwhile, we also discuss the nonlinear regression with $I(1)$
processes when $g(\cdot,\cdot)$ is asymptotically homogeneous. By using
Theorem~\ref{th3.2} in Section~\ref{sec3}, we obtain asymptotic
normality for the parametric estimator with a stochastic normalized
rate, which is comparable to Theorem~5.2 in PP. However, our derivation
is done under Markov perspective, which carries with it the potential
of extending the theory to nonlinear and nonstationary autoregressive
processes, which seems to be hard to do with the approach of PP.

The rest of this paper is organized as follows. Some preliminary
results about Markov theory (especially Harris recurrent Markov chain)
and function classes are introduced in Section~\ref{sec2}. The main
results of this paper and their extensions are given in Sections~\ref
{sec3} and \ref{sec4}, respectively. Some simulation studies are
carried out in Section~\ref{sec5} and the empirical application is
given in Section~\ref{sec6}. Section~\ref{sec7} concludes the paper.
The outline of the proofs of the main results is given in an \hyperref[app]{Appendix}.
The supplemental document [\citet{LTG15}] includes some
additional simulated examples, the detailed proofs of the main results
and the proofs of some auxiliary results.

\section{Preliminary results}\label{sec2}

To make the paper self-contained, in this section, we first provide
some basic definitions and preliminary results for a Harris recurrent
Markov process $\{X_t\}$, and then define function classes in a way
similar to those introduced in PP.

\subsection{Markov theory}\label{sec2.1}

Let $\{X_t, t\geq0\}$ be a $\phi$-irreducible Markov chain on the state
space $({\mathbb{E}}, {\mathcal{E}})$ with transition probability
$\mathsf{P}$. This means that for any set $A\in{\mathcal{E}}$ with
$\phi(A)>0$, we
have $\sum_{t=1}^\infty\mathsf{P}^t(x,A)>0$ for  $x\in {\mathbb{E}}$.
We further assume that the $\phi$-irreducible Markov chain $\{X_t\}$ is
Harris recurrent.

\begin{defn}\label{defn21}
A Markov chain $\{X_t\}$ is Harris
recurrent if,
for any set $B \in \varepsilon^+$ and given $X_0 = x$
for all $x \in \mathbb{E}$,
$\{X_t\}$ returns to $B$
infinitely often with probability one, where $\varepsilon^+$ is
defined as in \citet{KT01}.
\end{defn}

The Harris recurrence allows one to construct a split chain, which
decomposes the partial sum of functions of $\{X_t\}$ into blocks of
independent and identically distributed (i.i.d.) parts and two
asymptotically negligible remaining parts. Let $\tau_k$ be the
regeneration times, $n$ the number of observations and $N(n)$ the
number of regenerations as in \citet{KT01}, where
they use the notation $T(n)$ instead of $N(n)$. For the process $\{
G(X_t): t\geq0\}$, defining
\[
Z_{k}=\cases{ %
\displaystyle \sum
_{t=0}^{\tau_{0}}G(X_{t}), & \quad $k=0$,
\vspace*{3pt}
\cr
\displaystyle\sum_{t=\tau_{k-1}+1}^{\tau_{k}}G(X_{t}),
& \quad $1\leq k\leq N(n)$,\vspace*{3pt}
\cr
\displaystyle\sum
_{t=\tau_{N(n)}+1}^{n}G(X_{t}),&\quad $k=N(n)+1$,}
\]
where $G(\cdot)$ is a real function defined on $\mathbb{R}$, then
we have
%
\begin{equation}
\label{eq2.1} S_{n}(G)=\sum_{t=0}^nG(X_t)=Z_{0}+
\sum_{k=1}^{N(n)}Z_{k}+Z_{N(n)+1}.
\end{equation}
From \citet{N84}, we know that $\{Z_{k}, k\geq1\}$ is a sequence of
i.i.d. random variables, and $Z_{0}$ and $Z_{N(n)+1}$ converge to zero
almost surely (a.s.) when they are divided by the number of
regenerations $N(n)$ [using Lemma~3.2 in \citet{KT01}].

The general Harris recurrence only yields stochastic rates of
convergence in asymptotic theory of the parametric and nonparametric
estimators (see, e.g., Theorems \ref{th3.1} and \ref{th3.2} below),
where distribution and size of the number of regenerations $N(n)$ have
no a priori known structure but fully depend on the underlying process
$\{X_t\}$. To obtain a specific rate of $N(n)$ in our asymptotic theory
for the null recurrent process, we next impose some restrictions on the
tail behavior of the distribution of the recurrence times of the Markov chain.

\begin{defn}\label{defn22}
A Markov chain $\{X_t\}$ is $\beta$-null recurrent if there exist a~small
nonnegative function $f$, an initial measure $\lambda$, a constant
$\beta\in(0,1)$, and a slowly varying function $L_f(\cdot)$ such that
%
\begin{equation}
\label{eq2.2} \mathsf{E}_{\lambda} \Biggl(\sum_{t=1}^n
f(X_t) \Biggr)\sim\frac
{1}{\Gamma(1+\beta
)}n^{\beta}L_f(n),
\end{equation}
where $\mathsf{E}_{\lambda}$ stands for the expectation with initial
distribution $\lambda$ and $\Gamma(1+\beta)$ is
the Gamma function with parameter $1+\beta$.
\end{defn}

The definition of a small function $f$ in the above definition can be
found in some existing literature [c.f., page~15 in \citet{N84}].
Assuming $\beta$-null recurrence restricts the tail behavior of the
recurrence time of the process to be a regularly varying function. In
fact, for all small functions $f$, by Lemma~3.1 in \citet{KT01}, we can find an $L_s(\cdot)$ such that (\ref{eq2.2})
holds for the $\beta$-null recurrent Markov chain with $L_f(\cdot
)=\pi
_s(f)L_s(\cdot)$, where $\pi_s(\cdot)$ is an invariant measure of the
Markov chain $\{X_t\}$, $\pi_s(f)=\int f(x)\pi_s(dx)$ and $s$ is the
small function in the minorization inequality (3.4) of \citet{KT01}. Letting $L_s(n)=L_f(n)/(\pi_s(f))$ and following the
argument in \citet{KT01}, we may show that the
regeneration number $N(n)$ of the $\beta$-null recurrent Markov chain
$\{X_t\}$ has the following asymptotic distribution:
%
\begin{equation}
\label{eq2.3} \frac{N(n)}{n^{\beta}L_s(n)}\stackrel{d}\longrightarrow M_{\beta}(1),
\end{equation}
where $M_{\beta}(t)$, $t\geq0$ is the Mittag--Leffler process with
parameter $\beta$ [c.f., \citet{K84}]. Since $N(n)<n$ a.s. for the
null recurrent case by (\ref{eq2.3}), the rates of convergence for the
nonparametric kernel estimators are slower than those for the
stationary time series case [c.f., \citet{KMT07}, \citet{GKLT14}]. However, this is not necessarily the case for the
parametric estimator in our model (\ref{eq1.1}). In Section~\ref{sec3}
below, we will show that our rate of convergence in the null recurrent
case is slower than that for the stationary time series for integrable
$g(\cdot,\cdot)$ and may be faster than that for the stationary time
series case for asymptotically homogeneous $g(\cdot,\cdot)$. In
addition, our rates of convergence also depend on the magnitude of
$\beta$, which measures the recurrence times of the Markov chain $\{X_t\}$.

\subsection{Examples of \texorpdfstring{$\beta$}{beta}-null recurrent Markov chains}\label{sec2.2}

For a stationary or positive recurrent process, $\beta=1$. We next give
several examples of $\beta$-null recurrent Markov chains with $0<\beta<1$.

\begin{exa}[($1/2$-null recurrent Markov chain)]\label{exa21}
(i) Let a random walk process be defined as
%
\begin{equation}
\label{eq2.4} X_t=X_{t-1}+ x_t, \qquad t=1,2,
\ldots, X_0=0,
\end{equation}
where $\{x_t\}$ is a sequence of i.i.d. random variables with $\mathsf{E}[x_1]=0$, $0<\mathsf{E}[x_1^2]<\infty$ and
$\mathsf{E}[|x_1|^4]<\infty$, and the distribution of $x_t$ is absolutely
continuous (with respect to the Lebesgue measure) with the density
function $f_0(\cdot)$ satisfying $\inf_{x\in{\mathbb{C}}_0}f_0(x)>0$
for all compact sets ${\mathbb{C}}_0$. Some existing papers including
\citet{KR54} have shown that $\{X_t\}$ defined by
(\ref{eq2.4}) is a $1/2$-null recurrent Markov chain.

(ii) Consider a parametric TAR model of the form:
%
\begin{equation}
\label{eq2.5}\hspace*{8pt} X_t=\alpha_1 X_{t-1}I
(X_{t-1}\in{\mathbb{S}} )+\alpha_2 X_{t-1}I
\bigl(X_{t-1}\in{\mathbb{S}}^c \bigr)+ x_t,
\qquad X_0=0,
\end{equation}
where ${\mathbb{S}}$ is a compact subset of ${\mathbb{R}}$, ${\mathbb
{S}}^c$ is the complement of ${\mathbb{S}}$, $\alpha_2=1$, $-\infty
<\alpha_1<\infty$, $\{x_t\}$ satisfies the corresponding conditions in
Example~\ref{exa21}(i) above. Recently, \citet{GTY13} have shown that
such a TAR process $\{X_t\}$ is a $1/2$-null recurrent Markov chain.
Furthermore, we may generalize the TAR model~(\ref{eq2.5}) to
\[
X_t=H(X_{t-1},{\bolds\zeta})I (X_{t-1}\in{
\mathbb{S}} )+X_{t-1}I \bigl(X_{t-1}\in{\mathbb{S}}^c
\bigr)+ x_t,
\]
where $X_0=0$, $\sup_{x\in{\mathbb{S}}}|H(x,{\bolds\zeta
})|<\infty$ and ${\bolds\zeta}$ is a parameter vector. According to \citet{TTG10},
the above autoregressive process is also
a $1/2$-null recurrent Markov chain.
\end{exa}

\begin{exa}[($\beta$-null recurrent Markov chain with
$\beta\neq1/2$)]\label{exa22}
Let $\{x_t\}$ be a sequence of i.i.d. random variables taking positive
values, and $\{X_t\}$ be defined as
\[
X_t=\cases{ %
X_{t-1}-1, & \quad
$X_{t-1}>1$,\vspace*{3pt}
\cr
x_t, &\quad $X_{t-1}
\in[0,1]$,}
\]
for $t\geq1$, and $X_0=C_0$ for some positive constant $C_0$.
\citet{MKT12} prove that $\{X_t\}$ is $\beta$-null recurrent if
and only if
\[
\mathsf{P} \bigl([x_1]>n\bigr)\sim n^{-\beta}l^{-1}(n),
\qquad 0<\beta<1,
\]
where $[ \cdot]$ is the integer function and $l(\cdot)$ is a slowly
varying positive function.
\end{exa}

From the above examples, the $\beta$-null recurrent Markov chain
framework is not restricted to linear processes [see Example~\ref{exa21}(ii)].
Furthermore, such a null recurrent class has the invariance property
that if $\{X_t\}$ is $\beta$-null recurrent, then for a one-to-one
transformation ${\mathcal{T}}(\cdot)$, $\{{\mathcal{T}}(X_t)\}$ is also
$\beta$-null recurrent [c.f., \citet{TTG10}]. Such
invariance property does not hold for the $I(1)$ processes. For other
examples of the $\beta$-null recurrent Markov chain, we refer to
Example~1 in \citet{S11}. For some general conditions on diffusion
processes to ensure the Harris recurrence is satisfied, we refer to
\citet{HL03} and \citet{BP09}.

\subsection{Function classes}\label{sec2.3}

Similar to Park and Phillips (\citeyear{PP99}, \citeyear{PP01}), we consider two classes of
parametric nonlinear functions: integrable functions and asymptotically
homogeneous functions, which include many commonly-used functions in
nonlinear regression. Let $\|{\mathbf{A}}\|=\sqrt{\sum_{i=1}^q\sum_{j=1}^qa_{ij}^2}$ for ${\mathbf{A}}=(a_{ij})_{q\times q}$, and $\|
{\mathbf{a}}\|$ be the Euclidean norm of vector ${\mathbf{a}}$. A
function $h(x): {\mathbb{R}}\rightarrow{\mathbb{R}}^d$ is $\pi
_s$-integrable if
\[
\int_{\mathbb{R}}\bigl\|h(x)\bigr\|\pi_s(dx)<\infty,
\]
where $\pi_s(\cdot)$ is the invariant measure of the Harris recurrent
Markov chain $\{X_t\}$. When $\pi_s(\cdot)$ is differentiable such that
$\pi_s(dx)=p_s(x)\,dx$, $h(x)$ is $\pi_s$-integrable if and only if
$h(x)p_s(x)$ is integrable, where $p_s(\cdot)$ is the invariant density
function for $\{X_t\}$. For the random walk case as in Example~\ref{exa21}(i),
the $\pi_s$-integrability reduces to the conventional integrability as
$\pi_s(dx)=dx$.

\begin{defn}\label{defn23}
A $d$-dimensional vector function
$h(x, {\bolds\theta})$ is said to be integrable on $\Theta$ if for
each ${\bolds\theta}\in\Theta$, $h(x,{\bolds\theta})$
is $\pi_s$-integrable and there exist a neighborhood $\mathbb{B}_{\bolds\theta}$ of $\bolds\theta$ and $M: \mathbb{R}\rightarrow
\mathbb{R}$
bounded and $\pi_s$-integrable such that $\|h(x,{\bolds\theta
}^\prime)-h(x,{\bolds\theta})\|\leq\|{\bolds\theta
}^\prime
-{\bolds\theta}\|M(x)$ for any ${\bolds\theta}^\prime\in
{\mathbb{B}}_{\bolds\theta}$.
\end{defn}

The above definition is comparable to Definition~3.3 in PP. However, in
our definition, we do not need condition (b) in Definition~3.3 of their
paper, which makes the integrable function family in this paper
slightly more general. We next introduce a class of asymptotically
homogeneous functions.

\begin{defn}\label{defn24}
For a $d$-dimensional vector function
$h(x, {\bolds\theta})$, let $h(\lambda x,\break {\bolds\theta
})=\kappa(\lambda,{\bolds\theta})H(x,{\bolds\theta
})+R(x,\lambda,{\bolds\theta})$, where $\kappa(\cdot,\cdot)$ is
nonzero. $h(\lambda x,{\bolds\theta})$ is said to be
asymptotically homogeneous on $\Theta$ if the following two conditions
are satisfied: (i)~$H(\cdot,{\bolds\theta})$ is locally bounded
uniformly for any ${\bolds\theta}\in\Theta$ and continuous with
respect to ${\bolds\theta}$; (ii) the remainder term
$R(x,\lambda
,{\bolds\theta})$ is of order smaller than $\kappa(\lambda
,{\bolds\theta})$ as $\lambda\rightarrow\infty$ for any
${\bolds\theta}\in\Theta$. As in PP, $\kappa(\cdot,\cdot)$
is the
asymptotic order of $h(\cdot,\cdot)$ and $H(\cdot,\cdot)$ is the limit
homogeneous function.
\end{defn}

The above definition is quite similar to that of an $H$-regular
function in PP except that the regularity condition (a) in
Definition~3.5 of PP is replaced by the local boundness condition (i) in
Definition~\ref{defn24}. Following Definition~3.4 in PP, as $R(x,\lambda
,{\bolds\theta})$ is of order smaller than $\kappa(\cdot,\cdot)$,
we have either
%
\begin{equation}
\label{eq2.6} R (x,\lambda,{\bolds\theta} )=a(\lambda ,{\bolds\theta
})A_R (x,{\bolds\theta} )
\end{equation}
or
%
\begin{equation}
\label{eq2.7} R (x,\lambda,{\bolds\theta} )=b(\lambda ,{\bolds\theta
})A_R (x,{\bolds\theta} )B_R(\lambda x,{\bolds
\theta}),
\end{equation}
where $a(\lambda,{\bolds\theta})=o(\kappa(\lambda
,{\bolds
\theta}))$, $b(\lambda,{\bolds\theta})=O(\kappa(\lambda
,{\bolds\theta}))$ as $\lambda\rightarrow\infty$,
$\sup_{{\bolds\theta}\in\Theta}A_R(\cdot,{\bolds\theta
})$ is
locally bounded, and $\sup_{{\bolds\theta}\in\Theta}B_R(\cdot
,{\bolds\theta})$ is bounded and vanishes at infinity. 

Note that the above two definitions can be similarly generalized to the
case that $h(\cdot,\cdot)$ is a $d\times d$ matrix of functions.
Details are omitted here to save space. Furthermore, when the process
$\{X_t\}$ is positive recurrent, an asymptotically homogeneous function
$h(x, {\bolds\theta})$ might be also integrable on $\Theta$ as
long as the density function of the process $p_s(x)$ is integrable and
decreases to zero sufficiently fast when $x$ diverges to infinity.

\section{Main results}\label{sec3}

In this section, we establish some asymptotic results for the
parametric estimators of ${\bolds\theta}_0$ when $g(\cdot,\cdot)$
and its derivatives belong to the two classes of functions introduced
in Section~\ref{sec2.3}.

\subsection{Integrable function on \texorpdfstring{$\Theta$}{Theta}}\label{sec3.1}

We first consider estimating model (\ref{eq1.1}) by the NLS approach,
which is also used by PP in the unit root framework. Define the loss
function by
%
\begin{equation}
\label{eq3.1} L_{n,g}({\bolds\theta})=\sum
_{t=1}^n \bigl(Y_t-g(X_t,{
\bolds \theta}) \bigr)^2.
\end{equation}
We can obtain the resulting estimator $\widehat{\bolds\theta}_n$
by minimizing $L_{n,g}({\bolds\theta})$ over ${\bolds
\theta
}\in\Theta$, that is,
%
\begin{equation}
\label{eq3.2} \widehat{\bolds\theta}_n=\mathop{\arg
\min}_{{\bolds\theta}\in
\Theta
}L_{n,g}({\bolds\theta}).
\end{equation}

For ${\bolds\theta}=(\theta_1,\ldots,\theta_d)^\tau$, let
\[
\dot{g}(x,{\bolds\theta})= \biggl(\frac{\partial
g(x,{\bolds\theta})}{\partial\theta_j} \biggr)_{d\times1},
\qquad \ddot {g}(x,{\bolds\theta})= \biggl(\frac{\partial^2 g(x,{\bolds
\theta
})}{\partial\theta_i \partial\theta_j} \biggr)_{d\times d}.
\]
Before deriving the asymptotic properties of $\widehat{\bolds
\theta
}_n$ when $g(\cdot,\cdot)$ and its derivatives are integrable on
$\Theta
$, we give some regularity conditions.
%

\begin{ass}\label{ass31}
(i) $\{X_t\}$ is a Harris recurrent
Markov chain with invariant measure $\pi_s(\cdot)$.

(ii) $\{e_t\}$ is a sequence of i.i.d. random variables with mean
zero and finite variance $\sigma^2$, and is independent of $\{X_t\}$.
\end{ass}

\begin{ass}\label{ass32}
(i) $g(x,{\bolds\theta})$ is
integrable on $\Theta$, and for all ${\bolds\theta}\neq
{\bolds
\theta}_0$, $\int [g(x,{\bolds\theta})-g(x,{\bolds
\theta
}_0) ]^2\pi_s(dx)>0$.

(ii) Both $\dot{g}(x,{\bolds\theta})$ and $\ddot
{g}(x,{\bolds\theta})$ are integrable on $\Theta$, and the matrix
\[
\ddot{L}({\bolds\theta}):=\int\dot{g}(x,{\bolds\theta })\dot
{g}^\tau(x,{\bolds\theta})\pi_s(dx)
\]
is positive definite when ${\bolds\theta}$ is in a neighborhood of
${\bolds\theta}_0$.
\end{ass}

\begin{rem}\label{rem31}
In Assumption~\ref{ass31}(i), $\{X_t\}$ is assumed to
be Harris recurrent, which includes both the positive and null
recurrent Markov chains. The i.i.d. restriction on $\{e_t\}$ in
Assumption~\ref{ass31}(ii) may be replaced by the condition that $\{e_t\}$ is
an irreducible, ergodic and strongly mixing process with mean zero and
certain restriction on the mixing coefficient and moment conditions
[c.f., Theorem~3.4 in \citet{KMT07}]. Hence, under some
mild conditions, $\{e_t\}$ can include the well-known AR and ARCH
processes as special examples. However, for this case, the techniques
used in the proofs of Theorems \ref{th3.1} and \ref{th3.2} below need
to be modified by noting that the compound process $\{X_t,e_t\}$ is
Harris recurrent. Furthermore, the homoskedasticity on the error term
can also be relaxed, and we may allow the existence of certain
heteroskedasticity structure, that is, $e_t=\sigma(X_t)\eta_t$, where
$\sigma^2(\cdot)$ is the conditional variance function and $\{\eta
_t\}$
satisfies Assumption~\ref{ass31}(ii) with a unit variance. However, the
property of the function $\sigma^2(\cdot)$ would affect the convergence
rates given in the following asymptotic results. For example, to ensure
the validity of Theorem~\ref{th3.1}, we need to further assume that
$\sigma^2(\cdot)$ is $\pi_s$-integrable, which indicates that $\Vert
\dot{g}(x,{\bolds\theta})\Vert^2\sigma^2(x)$ is integrable on
$\Theta$. As in the literature [c.f., \citet{KMT07}], we
need to assume the independence between $\{X_t\}$ and $\{\eta_t\}$.


Assumption~\ref{ass32} is quite standard and similar to the corresponding
conditions in PP. In particular, Assumption~\ref{ass32}(i) is a key condition
to derive the global consistency of the NLS estimator $\widehat
{\bolds\theta}_n$.
\end{rem}

We next give the asymptotic properties of $\widehat{\bolds\theta
}_n$. The following theorem is applicable for both stationary (positive
recurrent) and nonstationary (null recurrent) time series.

\begin{thm}\label{th3.1}
Let Assumptions \ref{ass31} and \ref{ass32} hold.
\begin{longlist}[(a)]
\item[(a)] The solution $\widehat{\bolds\theta}_n$ which minimizes the
loss function $L_{n,g}({\bolds\theta})$ over $\Theta$ is
consistent, that is,
%
\begin{equation}
\label{eq3.3} \widehat{\bolds\theta}_n-{\bolds\theta}_0=o_P(1).
\end{equation}

\item[(b)]  The estimator $\widehat{\bolds\theta}_n$ has an asymptotically
normal distribution of the form:
%
\begin{equation}
\label{eq3.4} \sqrt{N(n)} (\widehat{\bolds\theta}_n-{\bolds\theta
}_0 )\stackrel{d}\longrightarrow{\mathsf{N}} \bigl({
\mathbf0}_d, \sigma^2\ddot{L}^{-1}({\bolds
\theta}_0) \bigr),
\end{equation}
where $\mathbf{0}_d$ is a $d$-dimensional null vector.
\end{longlist}
\end{thm}

\begin{rem}
Theorem~\ref{th3.1} shows that $\widehat
{\bolds\theta}_n$ is asymptotically normal with a stochastic
convergence rate $\sqrt{N(n)}$ for both the stationary and
nonstationary cases. However, $N(n)$ is usually unobservable and its
specific rate depends on $\beta$ and $L_s(\cdot)$ if $\{X_t\}$ is
$\beta$-null recurrent (see Corollary~\ref{cor32} below). We next discuss how to link
$N(n)$ with a directly observable hitting time. Indeed, if ${\mathbb
{C}}\in{\mathcal{E}}$ and $I_{\mathbb{C}}$ has a $\phi$-positive
support, the number of times that the process visits ${\mathbb{C}}$ up
to the time $n$ is defined by $N_{\mathbb{C}}(n)=\sum_{t=1}^nI_{\mathbb{C}}(X_t)$. By Lemma~3.2 in \citet{KT01}, we have
%
\begin{equation}
\label{eq3.5} \frac{N_{\mathbb{C}}(n)}{N(n)}\longrightarrow\pi_s({\mathbb{C}})
\qquad \mbox{a.s.},
\end{equation}
if $\pi_s({\mathbb{C}})=\pi_s I_{\mathbb{C}}=\int_{\mathbb C}\pi
_s(dx)<\infty$. A possible estimator of $\beta$ is
%
\begin{equation}
\label{eq3.5*} \widehat{\beta}=\frac{\ln N_{\mathbb{C}}(n)}{\ln n},
\end{equation}
which is strongly consistent as shown by \citet{KT01}. However, it is usually of somewhat limited practical use due to
the slow convergence rate [c.f., Remark~3.7 of \citet{KT01}]. A simulated example is given in Appendix B of the supplemental
document  to discuss the finite sample performance of the estimation
method in (\ref{eq3.5*}).
\end{rem}

By (\ref{eq3.5}) and Theorem~\ref{th3.1}, we can obtain the following
corollary directly.

\begin{cor}\label{cor31}
Suppose that the conditions of
Theorem~\ref{th3.1} are satisfied, and let ${\mathbb{C}}\in{\mathcal{E}}$ such
that $I_{\mathbb{C}}$ has a $\phi$-positive support and $\pi
_s({\mathbb
C})<\infty$. Then the estimator $\widehat{\bolds\theta}_n$ has an
asymptotically normal distribution of the form:
%
\begin{equation}
\label{eq3.6} \sqrt{N_{\mathbb{C}}(n)} (\widehat{\bolds\theta
}_n-{\bolds \theta}_0 )\stackrel{d}\longrightarrow{
\mathsf{N}} \bigl({\mathbf 0}_d, \sigma^2
\ddot{L}_{\mathbb{C}}^{-1}({\bolds\theta }_0) \bigr),
\end{equation}
where $\ddot{L}_{\mathbb{C}}({\bolds\theta}_0)=\pi
_s^{-1}({\mathbb
{C}})\ddot{L}({\bolds\theta}_0)$.
\end{cor}

\begin{rem}\label{rem33}
In practice, we may choose ${\mathbb C}$ as a
compact set such that $\phi({\mathbb C})>0$ and $\pi_s({\mathbb
C})<\infty$. In the additional simulation study (Example~B.1) given in
the supplemental document, for two types of $1/2$-null recurrent Markov
processes, we choose ${\mathbb C}=[-A,A]$ with the positive constant
$A$ carefully chosen, which works well in our setting. If $\pi_s(\cdot
)$ has a continuous derivative function $p_s(\cdot)$, we can show that
\[
\ddot{L}_{\mathbb{C}}({\bolds\theta}_0)=\int\dot {g}(x,{\bolds
\theta}_0)\dot{g}^\tau(x,{\bolds\theta}_0)p_{\mathbb
{C}}(x)
\,dx \qquad \mbox{with } p_{\mathbb{C}}(x)=p_s(x)/
\pi_{s}({\mathbb{C}}).
\]
The density function $p_{\mathbb{C}}(x)$ can be estimated by the kernel
method. Then, replacing ${\bolds\theta}_0$ by the NLS estimated
value, we can obtain a consistent estimate for $\ddot{L}_{\mathbb
{C}}({\bolds\theta}_0)$. Note that $N_{\mathbb{C}}(n)$ is
observable and $\sigma^2$ can be estimated by calculating the variance
of the residuals $\widehat{e}_t=Y_t-g(X_t,\widehat{\bolds\theta
}_n)$. Hence, for inference purposes,\vspace*{1pt} one may not need to estimate
$\beta$ and $L_s(\cdot)$ when $\{X_t\}$ is $\beta$-null recurrent, as
$\ddot{L}_{\mathbb{C}}({\bolds\theta}_0)$, $\sigma^2$ and
$N_{\mathbb{C}}(n)$ in (\ref{eq3.6}) can be explicitly computed without
knowing any information about $\beta$ and $L_s(\cdot)$.
\end{rem}






From (\ref{eq3.4}) in Theorem~\ref{th3.1} and (\ref{eq2.3}) in
Section~\ref{sec2.1} above, we have the following corollary.

\begin{cor}\label{cor32}
Suppose that the conditions of
Theorem~\ref{th3.1} are satisfied. Furthermore, $\{X_t\}$ is a $\beta
$-null recurrent Markov chain with $0<\beta<1$. Then we have
%
\begin{equation}
\label{eq3.7} \widehat{\bolds\theta}_n-{\bolds\theta}_0=O_P
\biggl(\frac
{1}{\sqrt{n^{\beta}L_s(n)}} \biggr),
\end{equation}
where $L_s(n)$ is defined in Section~\ref{sec2.1}.
\end{cor}

\begin{rem}\label{rem34}
As $\beta<1$ and $L_s(n)$ is a slowly varying
positive function, for the integrable case, the rate of convergence of
$\widehat{\bolds\theta}_n$ is slower than $\sqrt{n}$, the rate of
convergence of the parametric NLS estimator in the stationary time
series case. Combining (\ref{eq2.3}) and Theorem~\ref{th3.1}, the
result (\ref{eq3.7}) can be strengthened to
%
\begin{equation}
\label{eq3.7*} \sqrt{n^{\beta}L_s(n)} (\widehat{\bolds
\theta}_n-{\bolds \theta }_0)\stackrel{d}\longrightarrow
\bigl[M_\beta^{-1}(1)\sigma^2\ddot
{L}^{-1}({\bolds\theta}_0) \bigr]^{1/2}\cdot
\mathsf{N}_d,
\end{equation}
where $\mathsf{N}_d$ is a $d$-dimensional normal distribution with mean
zero and covariance matrix being the identity matrix, which is
independent of $M_\beta(1)$. A similar result is also obtained by \citet{CW12}. Corollary~\ref{cor32} and (\ref{eq3.7*}) complement the
existing results on the rates of convergence of nonparametric
estimators in $\beta$-null recurrent Markov processes [c.f., \citet{KMT07}, \citet{GKLT14}]. For the random walk
case, which corresponds to $1/2$-null recurrent Markov chain, the rate
of convergence is $n^{1/4}$, which is similar to a result obtained by
PP for the processes that are of $I(1)$ type.
\end{rem}

\subsection{Asymptotically homogeneous function on \texorpdfstring{$\Theta$}{Theta}}\label{sec3.2}

We next establish an asymptotic theory for a parametric estimator of
${\bolds\theta}_0$ when $g(\cdot,\cdot)$ and its derivatives
belong to a class of asymptotically homogeneous functions. For a unit
root process $\{X_t\}$, PP establish the consistency and limit
distribution of the NLS estimator $\widehat{\bolds\theta}_n$ by
using the local time technique. Their method relies on the linear
framework of the unit root process, the functional limit theorem of the
partial sum process and the continuous mapping theorem. The Harris
recurrent Markov chain is a general process and allows for a possibly
nonlinear framework, however. In particular, the null recurrent Markov
chain can be seen as a nonlinear generalization of the linear unit root
process. Hence, the techniques used by PP for establishing the
asymptotic theory is not applicable in such a possibly nonlinear Markov
chain framework. Meanwhile, as mentioned in Section~\ref{sec1}, the
methods used to prove Theorem~\ref{th3.1} cannot be applied here
directly because the asymptotically homogeneous functions usually are
not bounded and integrable. This leads to the violation of the
conditions in the ergodic theorem when the process is null recurrent.
In fact, most of the existing limit theorems for the null recurrent
Markov process $h(X_t)$ [c.f., Chen (\citeyear{C99}, \citeyear{C00})] only consider the
case where $h(\cdot)$ is bounded and integrable. Hence, it is quite
challenging to extend Theorem~3.3 in PP to the case of general null
recurrent Markov chains and establish an asymptotic theory for the NLS
estimator for the case of asymptotically homogeneous functions.

To address the above concerns, we have to modify the NLS estimator
$\widehat{\bolds\theta}_n$. Let $M_n$ be a positive and increasing
sequence satisfying $M_n\rightarrow\infty$ as $n\rightarrow\infty$, but
is dominated by a certain polynomial rate. We define the modified loss
function by
%
\begin{equation}
\label{eq3.8} Q_{n,g}({\bolds\theta})=\sum
_{t=1}^n \bigl[Y_t-g(X_t,{
\bolds \theta}) \bigr]^2I \bigl(|X_t|\leq M_n \bigr).
\end{equation}
The modified NLS (MNLS) estimator $\overline{\bolds\theta}_n$ can
be obtained by minimizing $Q_{n,g}({\bolds\theta})$ over
${\bolds\theta}\in\Theta$,
%
\begin{equation}
\label{eq3.9} \overline{\bolds\theta}_n=\mathop{\arg
\min}_{{\bolds\theta}\in
\Theta
}Q_{n,g}({\bolds\theta}).
\end{equation}

The above truncation technique enables us to develop the limit theorems
for the parametric estimate $\overline{\bolds\theta}_n$ even when
the function $g(\cdot,\cdot)$ or its derivatives are unbounded. A
similar truncation idea is also used by \citet{L07} to estimate the
ARMA-GARCH model when the second moment may not exist [it is called as
the self-weighted method by \citet{L07}]. However, Assumption~2.1 in
\citet{L07} indicates the stationarity for the model. The Harris
recurrence considered in the paper is more general and includes both
the stationary and nonstationary cases. As $M_n\rightarrow\infty$, for
the integrable case discussed in Section~\ref{sec3.1}, we can easily
show that $\overline{\bolds\theta}_n$ has the same asymptotic
distribution as $\widehat{\bolds\theta}_n$ under some regularity
conditions. In Example~\ref{exa51} below, we compare the finite sample
performance of these two estimators, and find that they are quite
similar. Furthermore, when $\{X_t\}$ is positive recurrent, as
mentioned in the last paragraph of Section~\ref{sec2.3}, although the
asymptotically homogeneous $g(x,{\bolds\theta})$ and its
derivatives are unbounded and not integrable on $\Theta$, it may be
reasonable to assume that $g(x, {\bolds\theta})p_s(x)$ and its
derivatives (with respect to ${\bolds\theta}$) are integrable on
$\Theta$. In this case, Theorem~\ref{th3.1} and Corollary~\ref{cor31} in Section~\ref
{sec3.1} still hold for the estimation $\overline{\bolds\theta}_n$
and the role of $N(n)$ [or $N_{\mathbb C}(n)]$ is the same as that of
the sample size, which implies that the root-$n$ consistency in the
stationary time series case can be derived. Hence, we only consider the
null recurrent $\{X_t\}$ in the remaining subsection.

Let
\begin{eqnarray*}
B_i(1) &=& [i-1,i ), \qquad i=1,2,\ldots,[M_n], \qquad
B_{[M_n]+1}(1)= \bigl[[M_n],M_n \bigr],
\\
B_i(2) &=& [-i,-i+1 ), \qquad i=1,2,\ldots,[M_n],
\qquad \! B_{[M_n]+1}(2)= \bigl[-M_n,-[M_n] \bigr].
\end{eqnarray*}
It\vspace*{1.5pt} is easy to check that $B_i(k)$, $i=1,2,\ldots,[M_n]+1$, $k=1,2$, are
disjoint, and $[-M_n,M_n]=\bigcup_{k=1}^2\bigcup_{i=1}^{[M_n]+1}B_i(k)$. Define
\begin{eqnarray*}
&& \zeta_0(M_n)=\sum_{i=0}^{[M_n]}
\bigl[\pi_s\bigl(B_{i+1}(1)\bigr)+\pi _s
\bigl(B_{i+1}(2)\bigr) \bigr],
\end{eqnarray*}
and
\begin{eqnarray*}
{\bolds\Lambda}_n({\bolds\theta}) &=& \sum
_{i=0}^{[M_n]}\dot {h}_g\biggl(
\frac{i}{M_n},{\bolds\theta}\biggr)\dot{h}_g^\tau
\biggl(\frac
{i}{M_n},{\bolds\theta}\biggr)\pi_s
\bigl(B_{i+1}(1)\bigr)
\\
&&{}+\sum_{i=0}^{[M_n]}\dot{h}_g
\biggl(\frac{-i}{M_n},{\bolds\theta }\biggr)\dot {h}_g^\tau
\biggl(\frac{-i}{M_n},{\bolds\theta}\biggr)\pi _s
\bigl(B_{i+1}(2)\bigr),
\\
\widetilde{\Lambda}_n({\bolds\theta},{\bolds\theta}_0)
&=& \sum_{i=0}^{[M_n]} \biggl[h_g
\biggl(\frac{i}{M_n},{\bolds\theta }\biggr)-h_g\biggl(
\frac
{i}{M_n},{\bolds\theta}_0\biggr) \biggr]^2
\pi_s\bigl(B_{i+1}(1)\bigr)
\\
&&{}+\sum_{i=0}^{[M_n]} \biggl[h_g
\biggl(\frac{-i}{M_n},{\bolds\theta }\biggr)-h_g\biggl(
\frac{-i}{M_n},{\bolds\theta}_0\biggr) \biggr]^2\pi
_s\bigl(B_{i+1}(2)\bigr),
\end{eqnarray*}
where $h_g(\cdot,\cdot)$ and $\dot{h}_g(\cdot,\cdot)$ will be defined
in Assumption~\ref{ass33}(i) below.

Some additional assumptions are introduced below to establish
asymptotic properties for $\overline{\bolds\theta}_n$.

\begin{ass}\label{ass33}
(i) $g(x,{\bolds\theta})$, $\dot
{g}(x,{\bolds\theta})$ and $\ddot{g}(x,{\bolds\theta})$ are
asymptotically homogeneous on $\Theta$ with asymptotic orders $\kappa
_g(\cdot)$, $\dot{\kappa}_g(\cdot)$ and $\ddot{\kappa}_g(\cdot
)$, and
limit homogeneous functions $h_g(\cdot,\cdot)$, $\dot{h}_g(\cdot
,\cdot
)$ and $\ddot{h}_g(\cdot,\cdot)$, respectively. Furthermore, the
asymptotic orders $\kappa_g(\cdot)$, $\dot{\kappa}_g(\cdot)$ and
$\ddot
{\kappa}_g(\cdot)$ are independent of ${\bolds\theta}$.

(ii) The function $h_g(\cdot, {\bolds\theta})$ is
continuous on the interval $[-1,1]$ for all ${\bolds\theta}\in
\Theta$. For all ${\bolds\theta}\neq{\bolds\theta}_0$, there
exist a continuous $\widetilde{\Lambda}(\cdot,{\bolds\theta}_0)$
which achieves unique minimum at ${\bolds\theta}={\bolds
\theta
}_0$ and a sequence of positive numbers $\{\widetilde{\zeta}(M_n)\}$
such that
%
\begin{equation}
\label{eq3.9*} \lim_{n\rightarrow\infty}\frac{1}{\widetilde{\zeta
}(M_n)}\widetilde {
\Lambda}_n({\bolds\theta},{\bolds\theta}_0)=\widetilde {
\Lambda}({\bolds\theta},{\bolds\theta}_0).
\end{equation}
For ${\bolds\theta}$ in a neighborhood of ${\bolds\theta}_0$,
both $\dot{h}_g(\cdot,{\bolds\theta})$ and $\ddot{h}_g(\cdot
,{\bolds\theta})$ are continuous on the interval $[-1,1]$ and
there exist a continuous and positive definite matrix ${\bolds\Lambda
}({\bolds\theta})$ and a sequence of positive numbers $\{\zeta
(M_n)\}$ such that
%
\begin{equation}
\label{eq3.9**} \lim_{n\rightarrow\infty}\frac{1}{\zeta(M_n)}{\bolds\Lambda
}_n({\bolds\theta})={\bolds\Lambda}({\bolds\theta}).
\end{equation}
Furthermore, both $\zeta_0(M_n)/\widetilde{\zeta}(M_n)$ and $\zeta
_0(M_n)/\zeta(M_n)$ are bounded, and $\zeta_0(l_n)/\zeta_0(M_n)=o(1)$
for $l_n\rightarrow\infty$ but $l_n=o(M_n)$.

(iii) The asymptotic orders $\kappa_g(\cdot)$, $\dot{\kappa
}_g(\cdot)$ and $\ddot{\kappa}_g(\cdot)$ are positive and nondecreasing
such that $\kappa_g(n)+\ddot{\kappa}_g(n)=O (\dot{\kappa
}_g(n) )$
as $n\rightarrow\infty$.

(iv) For each $x\in[-M_n,M_n]$, ${\mathcal N}_x(1):= \{y:
x-1<y<x+1 \}$ is a small set and the invariant density function
$p_s(x)$ is bounded away from zero and infinity.
\end{ass}

\begin{rem}\label{rem35}
Assumption~\ref{ass33}(i) is quite standard; see, for
example, condition (b) in Theorem~5.2 in PP. The restriction that the
asymptotic orders are independent of ${\bolds\theta}$ can be
relaxed at the cost of more complicated assumptions and more lengthy
proofs. For example, to ensure the global consistency of $\overline
{\bolds\theta}_n$, we need to assume that there exist $\epsilon
_\ast>0$ and a neighborhood ${\mathbb B}_{{\bolds\theta}_1}$ of
${\bolds\theta}_1$ for any ${\bolds\theta}_1\neq
{\bolds
\theta}_0$ such that
\[
\inf_{|p-\bar{p}|<\epsilon_\ast, |q-\bar{q}|<\epsilon_\ast} \inf_{{\bolds\theta}\in{\mathbb B}_{{\bolds\theta}_1}} \bigl|p\kappa
_g(n,{\bolds\theta})-q\kappa_g(n,{\bolds
\theta}_0) \bigr|\rightarrow\infty
\]
for $\bar{p},\bar{q}>0$. And to establish the asymptotic normality of
$\overline{\bolds\theta}_n$, we need to impose additional
technical conditions on the asymptotic orders and limit homogeneous
functions, similar to condition (b) in Theorem~5.3 of PP. The explicit
forms of ${\bolds\Lambda}({\bolds\theta})$, $\zeta(M_n)$,
$\widetilde{\Lambda}({\bolds\theta},{\bolds\theta}_0)$ and
$\widetilde{\zeta}(M_n)$ in Assumption~\ref{ass33}(ii) can be derived for some
special cases. For example, when $\{X_t\}$ is generated by a random
walk process, we have $\pi_s(dx)=dx$ and
\begin{eqnarray}
{\bolds\Lambda}_n({\bolds\theta})&=&\bigl(1+o(1)\bigr)\sum
_{i=-[M_n]}^{[M_n]}\dot{h}_g\biggl(
\frac{i}{M_n},{\bolds\theta}\biggr)\dot {h}_g^\tau
\biggl(\frac{i}{M_n},{\bolds\theta}\biggr)
\nonumber
\\
&=&\bigl(1+o(1)\bigr)M_n\int_{-1}^1
\dot{h}_g(x,{\bolds\theta})\dot {h}_g^\tau
(x,{\bolds\theta})\,dx,
\nonumber
\end{eqnarray}
which implies that $\zeta(M_n)=M_n$ and ${\bolds\Lambda
}({\bolds
\theta})=\int_{-1}^1\dot{h}_g(x,{\bolds\theta})\dot
{h}_g^\tau
(x,{\bolds\theta})\,dx$ in (\ref{eq3.9**}). The explicit forms of
$\widetilde{\Lambda}({\bolds\theta},{\bolds\theta}_0)$ and
$\widetilde{\zeta}(M_n)$ can be derived similarly for the above two
cases and details are thus omitted.
\end{rem}

Define ${\mathbf{J}}_g(n,{\bolds\theta}_0)=\dot{\kappa
}_g^2(M_n)\zeta(M_n){\bolds\Lambda}({\bolds\theta}_0)$. We next
establish an asymptotic theory for $\overline{\bolds\theta}_n$
when $\{X_t\}$ is null recurrent.

\begin{thm}\label{th3.2}
Let $\{X_t\}$ be a null recurrent Markov process, Assumptions~\ref{ass31}\textup{(ii)}
and \ref{ass33} hold.
\begin{longlist}[(a)]
\item[(a)] The solution $\overline{\bolds\theta}_n$ which minimizes the
loss function $Q_{n,g}({\bolds\theta})$ over $\Theta$ is
consistent, that is,
%
\begin{equation}
\label{eq3.10} \overline{\bolds\theta}_n-{\bolds\theta}_0=o_P(1).
\end{equation}

\item[(b)]  The estimator $\overline{\bolds\theta}_n$ has the
asymptotically normal distribution,
%
\begin{equation}
\label{eq3.11} N^{1/2}(n){\mathbf{J}}_g^{1/2}(n,{
\bolds\theta}_0) (\overline {\bolds\theta}_n-{\bolds
\theta}_0 )\stackrel {d}\longrightarrow{\mathsf{N}} \bigl({
\mathbf0}_d, \sigma ^2I_d \bigr),
\end{equation}
where $I_d$ is a $d\times d$ identity matrix.
\end{longlist}
\end{thm}

\begin{rem}\label{rem36}
From Theorem~\ref{th3.2}, the asymptotic
distribution of $\overline{\bolds\theta}_n$ for the asymptotically
homogeneous regression function is quite different from that of
$\widehat{\bolds\theta}_n$ for the integrable regression function
when the process is null recurrent. Such finding is comparable to those
in PP. The choice of $M_n$ in the estimation method and asymptotic
theory will be discussed in Corollaries \ref{cor33} and \ref{cor34} below.
\end{rem}

\begin{rem}\label{rem37}
As in Corollary~\ref{cor31}, we can modify (\ref
{eq3.11}) for inference purposes. Define ${\mathbf{J}}_{g,{\mathbb
{C}}}(n,{\bolds\theta}_0)=\dot{\kappa}_g^2(M_n){\bolds\Lambda
}_{n,\mathbb{C}}({\bolds\theta}_0)$, where
\begin{eqnarray*}
{\bolds\Lambda}_{n,{\mathbb{C}}}({\bolds\theta})&=&\sum
_{i=0}^{[M_n]}\dot{h}_g\biggl(
\frac{i}{M_n},{\bolds\theta}\biggr)\dot {h}_g^\tau
\biggl(\frac{i}{M_n},{\bolds\theta}\biggr)\frac{\pi_s(B_{i+1}(1))}{\pi
_s({\mathbb{C}})}
\\
&&{}+\sum_{i=0}^{[M_n]}\dot{h}_g
\biggl(\frac{-i}{M_n},{\bolds\theta }\biggr)\dot {h}_g^\tau
\biggl(\frac{-i}{M_n},{\bolds\theta}\biggr)\frac{\pi
_s(B_{i+1}(2))}{\pi_s({\mathbb{C}})},
\end{eqnarray*}
where ${\mathbb C}$ satisfies the conditions in Corollary~\ref{cor31}. Then, by
(\ref{eq3.5}) and (\ref{eq3.11}), we can show that
%
\begin{equation}
\label{eq3.12} N_{\mathbb{C}}^{1/2}(n){\mathbf{J}}_{g,{\mathbb
{C}}}^{1/2}(n,{
\bolds\theta}_0) (\overline{\bolds \theta }_n-{\bolds
\theta}_0 )\stackrel{d}\longrightarrow {\mathsf {N}} \bigl({
\mathbf0}_d, \sigma^2I_d \bigr).
\end{equation}
When $\{X_t\}$ is $\beta$-null recurrent, we can use the asymptotically
normal distribution theory (\ref{eq3.12}) to conduct statistical
inference without knowing any information of $\beta$ as $N_{\mathbb
{C}}(n)$ is observable and ${\mathbf{J}}_{g,{\mathbb
{C}}}(n,{\bolds
\theta}_0)$ can be explicitly computed through replacing ${\bolds
\Lambda}_{n,{\mathbb{C}}}({\bolds\theta}_0)$ by the plug-in
estimated value.
\end{rem}





From (\ref{eq3.11}) in Theorem~\ref{th3.2} and (\ref{eq2.3}) in
Section~\ref{sec2} above, we have the following two corollaries. The
rate of
convergence in (\ref{eq3.13}) below is quite general for $\beta$-null
recurrent Markov processes. When $\beta=1/2$, it is the same as the
convergence rate in Theorem~5.2 of PP.

\begin{cor}\label{cor33}
Suppose that the conditions of
Theorem~\ref{th3.2} are satisfied. Furthermore, let $\{X_t\}$ be a
$\beta$-null recurrent Markov chain with $0<\beta<1$. Taking
$M_n=M_0n^{1-\beta}L_s^{-1}(n)$ for some positive constant $M_0$, we have
%
\begin{equation}
\label{eq3.13} \overline{\bolds\theta}_n-{\bolds\theta}_0=O_P
\bigl(\bigl(n\dot {\kappa}_g^2(M_n)
\bigr)^{-1/2} \bigr).
\end{equation}
\end{cor}

\begin{cor}\label{cor34}
Suppose that the conditions of
Theorem~\protect\ref{th3.2} are satisfied. Let $g(x,{\bolds\theta
}_0)=x\theta_0$, $\{X_t\}$ be a random walk process and
$M_n=M_0n^{1/2}$ for some positive constant $M_0$. Then we have
%
\begin{equation}
\label{eq3.14} \overline{\theta}_{n}-\theta_0=O_P
\bigl(n^{-1} \bigr),
\end{equation}
where $\overline{\theta}_{n}$ is the MNLS estimator of $\theta_0$.
Furthermore,
\begin{equation}
\label{eq3.15} \sqrt{M_0^3N(n)n^{3/2}}
(\overline{\theta}_{n}-\theta_0 )\stackrel{d}
\longrightarrow{\mathsf{N}} \bigl(0,3\sigma^2/2 \bigr).
\end{equation}
\end{cor}

\begin{rem}\label{rem38}
For the simple linear regression model with
regressors generated by a random walk process, (\ref{eq3.14}) and
(\ref
{eq3.15}) imply the existence of super consistency. Corollaries \ref{cor33} and
\ref{cor34} show that the rates of convergence for the parametric estimator in
nonlinear cointegrating models rely not only on the properties of the
function $g(\cdot,\cdot)$, but also on the magnitude of $\beta$.

In the above two corollaries, we give the choice of $M_n$ for some
special cases. In fact, for the random walk process $\{X_t\}$ defined
as in Example~\ref{exa21}(i) with ${\sf E}[x_1^2]=1$, we have
\[
\frac{1}{\sqrt{n}}X_{[nr]}=\frac{1}{\sqrt{n}}\sum
_{i=1}^{[nr]} x_i\Rightarrow B(r),
\]
where $B(r)$ is a standard Brownian motion and ``$\Rightarrow$'' denotes
the weak convergence. Furthermore, by the continuous mapping theorem
[c.f., \citet{B68}],
\[
\sup_{0\leq r\leq1}\frac{1}{\sqrt{n}}X_{[nr]}\Rightarrow\sup
_{0\leq
r\leq1}B(r),
\]
which implies that it is reasonable to let $M_n=C_\alpha n^{1/2}$,
where $C_{\alpha}$ may be chosen such that
%
\begin{equation}
\label{eq3.16} \alpha= {\mathsf{P}} \Bigl(\sup_{0\leq r\leq1}B(r)\geq
C_{\alpha
} \Bigr)={\mathsf{P}} \bigl(\bigl\llvert B(1)\bigr\rrvert \geq
C_{\alpha} \bigr) = 2 \bigl(1- \Phi (C_{\alpha} ) \bigr),
\end{equation}
where the second equality is due to the reflection principle and $\Phi
(x) = \int_{-\infty}^x  (e^{-{u^2}/{2}}/\sqrt{2\pi} )\, du$.\vspace*{1pt}
This implies that $C_{\alpha}$ can be obtained when $\alpha$ is given,
such as $\alpha= 0.05$. For the general $\beta$-null recurrent Markov
process, the choice of the optimal $M_n$ remains as an open problem. We
conjecture that it may be an option to take $M_n=\widetilde
{M}n^{1-\widehat{\beta}}$ with $\widehat{\beta}$ defined in (\ref
{eq3.5*}) and $\widetilde{M}$ chosen by a data-driven method, and will
further study this issue in future research.
\end{rem}

\section{Discussions and extensions}\label{sec4}

In this section, we discuss the applications of our asymptotic results
in estimating the nonlinear heteroskedastic regression and nonlinear
regression with $I(1)$ processes. Furthermore, we also discuss possible
extensions of our model to the cases of multivariate regressors and
nonlinear autoregression.

\subsection{Nonlinear heteroskedastic regression}\label{sec4.1}

$\!$We introduce an estimation method for a parameter vector involved in
the conditional variance function. For simplicity, we consider the
model defined by
%
\begin{equation}
\label{eq4.1} Y_t=\sigma(X_t,{\bolds
\gamma}_0)e_{t*} \qquad \mbox{for } {\bolds
\gamma}_0\in\Upsilon\subset{\mathbb{R}}^{p},
\end{equation}
where $\{e_{t*}\}$ satisfies Assumption~\ref{ass31}(ii) with a unit variance,
$\sigma^2(\cdot,\cdot):{\mathbb{R}}^{p+1}\rightarrow{\mathbb{R}}$ is
positive, and ${\bolds\gamma}_0$ is the true value of the
$p$-dimensional parameter vector involved in the conditional variance
function. Estimation of the parametric nonlinear variance function
defined in (\ref{eq4.1}) is important in empirical applications as many
scientific studies depend on understanding the variability of the data.
When the covariates are integrated, \citet{HP12} study the
maximum likelihood estimation of the parameters in the ARCH and GARCH
models. A recent paper by \citet{HK13} further considers
the quasi maximum likelihood estimation in the GARCH-X models with
stationary and nonstationary covariates. We next consider the general
Harris recurrent Markov process $\{X_t\}$ and use a robust estimation
method for model (\ref{eq4.1}).

Letting $\varpi_0$ be a positive number such that ${\mathsf{E}}
[\log(e_{t*}^2) ]=\log(\varpi_0)$,
we have
%
\begin{eqnarray}
\log\bigl(Y_t^2\bigr)&=&\log \bigl(\sigma^2(X_t,{
\bolds\gamma}_0) \bigr)+\log\bigl(e_{t*}^2
\bigr)
\nonumber
\\
\label{eq4.2} &=&\log \bigl(\sigma^2(X_t,{\bolds
\gamma}_0) \bigr)+\log (\varpi _0)+\log
\bigl(e_{t*}^2\bigr)-\log(\varpi_0)
\\
&=:&\log \bigl(\varpi_0\sigma^2(X_t,{\bolds
\gamma}_0) \bigr)+\zeta _t,
\nonumber
\end{eqnarray}
where ${\mathsf{E}}(\zeta_t)=0$. Since our main interest lies in the
discussion of the asymptotic theory for the estimator of ${\bolds
\gamma}_0$, we first assume that $\varpi_0$ is known to simplify our
discussion. Model (\ref{eq4.2}) can be seen as another nonlinear mean
regression model with parameter vector ${\bolds\gamma}_0$ to be
estimated. The log-transformation would make data less skewed, and thus
the resulting volatility estimator may be more robust in terms of
dealing with heavy-tailed $\{e_{t*}\}$. Such transformation has been
commonly used to estimate the variability of the data in the stationary
time series case [c.f., \citet{PY03}, \citet{G07}, \citet{CCP09}]. However, any extension to Harris recurrent Markov chains
which may be nonstationary has not been done in the literature.

Our estimation method will be constructed based on (\ref{eq4.2}).
Noting that $\varpi_0$ is assumed to be known, define
%
\begin{equation}
\label{eq4.3} \hspace*{6pt}\sigma_*^2(X_t,{\bolds
\gamma}_0)=\varpi_0\sigma ^2(X_t,{
\bolds \gamma}_0) \quad \mbox{and} \quad g_*(X_t,{\bolds
\gamma}_0)=\log \bigl(\sigma_*^2(X_t,{\bolds
\gamma}_0)\bigr).
\end{equation}

\begin{longlist}[\textit{Case} (II).]
\item[\textit{Case} (I).] If $g_*(X_t,{\bolds\gamma})$ and its derivatives
are integrable on $\Upsilon$, the log-transformed nonlinear least
squares (LNLS) estimator $\widehat{\bolds\gamma}_n$ can be
obtained by minimizing $L_{n,\sigma}({\bolds\gamma})$ over
${\bolds\gamma}\in\Upsilon$, where
%
\begin{equation}
\label{eq4.4} L_{n,\sigma}({\bolds\gamma})=\sum
_{t=1}^n \bigl[\log \bigl(Y_t^2
\bigr)-g_*(X_t,{\bolds\gamma}) \bigr]^2.
\end{equation}
Letting Assumptions \ref{ass31} and \ref{ass32} be satisfied with $e_t$ and $g(\cdot
,\cdot)$ replaced by $\zeta_t$ and $g_*(\cdot,\cdot)$, respectively,
then the asymptotic results developed in Section~\ref{sec3.1} still
hold for $\widehat{\bolds\gamma}_n$.

\item[\textit{Case} (II).] If $g_*(X_t,{\bolds\gamma})$ and its derivatives
are asymptotically homogeneous on~$\Upsilon$, the log-transformed
modified nonlinear least squares (LMNLS) estimator $\overline
{\bolds\gamma}_n$ can be obtained by minimizing $Q_{n,\sigma
}({\bolds\gamma})$ over ${\bolds\gamma}\in\Upsilon$, where
%
\begin{equation}
\label{eq4.5} Q_{n,\sigma}({\bolds\gamma})=\sum
_{t=1}^n \bigl[\log \bigl(Y_t^2
\bigr)-g_\ast (X_t,{\bolds\gamma}) \bigr]^2I
\bigl(|X_t|\leq M_n \bigr),
\end{equation}
where $M_n$ is defined as in Section~\ref{sec3.2}. Then the asymptotic
results developed in Section~\ref{sec3.2} still hold for $\overline
{\bolds\gamma}_n$ under some regularity conditions such as a
slightly modified version of Assumptions \ref{ass31} and \ref{ass33}. Hence, \textit
{it is
possible to achieve the super-consistency result for $\overline
{\bolds\gamma}_n$ when $\{X_t\}$ is null recurrent}. 
\end{longlist}

In practice, however, $\varpi_0$ is usually unknown and needs to be
estimated. We next briefly discuss this issue for case (ii). We may
define the loss function by
\[
Q_{n}({\bolds\gamma},\varpi)=\sum_{t=1}^n
\bigl[\log \bigl(Y_t^2\bigr)-\log \bigl(\varpi
\sigma^2(X_t,{\bolds\gamma}) \bigr) \bigr]^2I
\bigl(|X_t|\leq M_n \bigr).
\]
Then the estimators $\overline{\bolds\gamma}_n$ and $\overline
{\varpi}_n$ can be obtained by minimizing $Q_{n}({\bolds\gamma
},\varpi)$ over ${\bolds\gamma}\in\Upsilon$ and $\varpi\in
{\mathbb
{R}}^+$. A simulated example (Example~B.2) is given in Appendix~B of
the supplemental document  to examine the finite sample performance of
the LNLS and LMNLS estimations considered in cases (i) and (ii), respectively.

\subsection{Nonlinear regression with $I(1)$ processes}\label{sec4.2}

As mentioned before, PP consider the nonlinear regression (\ref{eq1.1})
with the regressors $\{X_t\}$ generated by
%
\begin{equation}
\label{eq4.6} X_t=X_{t-1}+x_t, \qquad
x_t=\sum_{j=0}^{\infty}
\phi_j\varepsilon_{t-j},
\end{equation}
where $\{\varepsilon_j\}$ is a sequence of i.i.d. random variables and
$\{\phi_j\}$ satisfies some summability conditions. For simplicity, we
assume that $X_0=0$ throughout this subsection. PP establish a suite of
asymptotic results for the NLS estimator of the parameter ${\bolds
\theta}_0$ involved in (\ref{eq1.1}) when $\{X_t\}$ is defined by
(\ref
{eq4.6}). An open problem is how to establish such results by using the
$\beta$-null recurrent Markov chain framework. This is quite
challenging as $\{X_t\}$ defined by (\ref{eq4.6}) is no longer a Markov
process except for some special cases (for example, $\phi_j=0$ for
$j\geq1$).

We next consider solving this open problem for the case where $g(\cdot
,\cdot)$ is asymptotically homogeneous on $\Theta$ and derive an
asymptotic theory for $\overline{\bolds\theta}_n$ by using
Theorem~\ref{th3.2} (the discussion for the integrable case is more
complicated, and will be considered in a future study). Our main idea
is to approximate $X_t$ by $X_t^\ast$ which is defined by
\[
X_t^*=\phi\sum_{s=1}^t
\varepsilon_{s}, \qquad \phi:=\sum_{j=0}^\infty
\phi _j\neq0,
\]
and then show that the asymptotically homogeneous function of $X_t$ is
asymptotically equivalent to the same function of $X_t^\ast$. As $\{
X_t^\ast\}$ is a random walk process under the Assumption~E.1 (see
Appendix~E of the supplemental document), we can then make use of
Theorem~\ref{th3.2}. Define
%
\begin{equation}
\label{eq4.7} {\mathbf{J}}_{g*}(n,{\bolds\theta}_0)=
\dot{\kappa }_g^2(M_n)M_n \biggl(
\int_{-1}^1 \dot{h}_g(x,{\bolds
\theta}_0)\dot{h}_g^\tau(x,{\bolds \theta
}_0)\,dx \biggr).
\end{equation}
We next give some asymptotic results for $\overline{\bolds\theta
}_n$ for the case where $\{X_t\}$ is a unit root process (\ref{eq4.6}),
and the proof is provided in Appendix E of the supplemental document.

\begin{thm}\label{th4.1}
Let\vspace*{1pt} Assumptions \textup{E.1} and \textup{E.2} in Appendix \textup{E} of the supplemental document
hold, and $n^{-{1}/({2(2+\delta)})}M_n\rightarrow\infty$, where
$\delta>0$ is defined in Assumption \textup{E.1(i)}.
\begin{longlist}[(a)]
\item[(a)]  The solution $\overline{\bolds\theta}_n$ which minimizes the
loss function $Q_{n,g}({\bolds\theta})$ over $\Theta$ is
consistent, that is,
%
\begin{equation}
\label{eq4.8} \overline{\bolds\theta}_n-{\bolds\theta}_0=o_P(1).
\end{equation}

\item[(b)]  The estimator $\overline{\bolds\theta}_n$ has the
asymptotically normal distribution,
%
\begin{equation}
\label{eq4.9} N_\varepsilon^{1/2}(n){\mathbf{J}}_{g*}^{1/2}(n,{
\bolds\theta }_0) (\overline{\bolds\theta}_n-{\bolds
\theta }_0 )\stackrel{d}\longrightarrow{\mathsf{N}} \bigl({
\mathbf0}_d, \sigma ^2I_d \bigr),
\end{equation}
where $N_\varepsilon(n)$ is the number of regenerations for the random
walk $\{X_t^*\}$.
\end{longlist}
\end{thm}

\begin{rem}\label{rem41}
Theorem~\ref{th4.1} establishes an asymptotic
theory for $\overline{\bolds\theta}_n$ when $\{X_t\}$ is a unit
root process (\ref{eq4.6}). Our results are comparable with Theorems
5.2 and 5.3 in PP. However, we establish asymptotic normality in (\ref
{eq4.9}) with stochastic rate $N_\varepsilon^{1/2}(n){\mathbf
{J}}_{g*}^{1/2}(n,{\bolds\theta}_0)$, and PP establish their
asymptotic mixed normal distribution theory with a deterministic rate.
As $N_\varepsilon^{1/2}(n){\mathbf{J}}_{g*}^{1/2}(n,{\bolds
\theta
}_0)\propto n^{1/4}{\mathbf{J}}_{g*}^{1/2}(n,{\bolds\theta}_0)$ in
probability, if we take $M_n=M_0\sqrt{n}$ as in Corollary~\ref{cor34}, we will
find that our rate of convergence of $\overline{\bolds\theta}_n$
is the same as that derived by PP.
\end{rem}

\subsection{Extensions to multivariate regression and nonlinear
autoregression}\label{sec4.3}

The theoretical results developed in Section~\ref{sec3} are limited to
nonlinear regression with a univariate Markov process. A natural
question is whether it is possible to extend them to the more general
case with multivariate covariates. In the unit root framework, it is
well known that it is difficult to derive the limit theory for the case
of multivariate unit root processes, as the vector Brownian motion is
transient when the dimension is larger than (or equal to) $3$. In
contrast, under the framework of the Harris recurrent Markov chains, it
is possible for us to generalize the theoretical theory to the
multivariate case (with certain restrictions). For example, it is
possible to extend the theoretical results to the case with one
nonstationary regressor and several other stationary regressors. We
next give an example of vector autoregressive (VAR) process which may
be included in our framework under certain conditions.

\begin{exa}\label{exa41}
Consider a $q$-dimensional VAR(1) process $\{
{\mathbf X}_t\}$ which is defined by
%
\begin{equation}
\label{eq4.10} {\mathbf X}_t={\mathbf A} {\mathbf X}_{t-1}+{
\mathbf b}+{\mathbf x}_t, \qquad t=1,2,\ldots,
\end{equation}
where ${\mathbf X}_0={\mathbf0}_q$, ${\mathbf A}$ is a $q\times q$ matrix,
${\mathbf b}$ is a $q$-dimensional vector and $\{{\mathbf x}_t\}$ is a
sequence of i.i.d. $q$-dimensional random vectors with mean zero. If
all the eigenvalues of the matrix ${\mathbf A}$ are inside the unit
circle, under some mild conditions on $\{{\mathbf x}_t\}$, Theorem~3 in
\citet{MKT12} shows that the VAR(1) process $\{
{\mathbf X}_t\}$ in (\ref{eq4.10}) is geometric ergodic, which belongs
to the category of positive recurrence. On the other hand, if the
matrix ${\mathbf A}$ has exactly one eigenvalue on the unit circle,
under some mild conditions on $\{{\mathbf x}_t\}$ and ${\mathbf b}$,
Theorem~4 in \citet{MKT12} shows that the VAR(1)
process $\{{\mathbf X}_t\}$ in (\ref{eq4.10}) is $\beta$-null recurrent
with $\beta=1/2$. For this case, the asymptotic theory developed in
Section~\ref{sec3} is applicable. However, when ${\mathbf A}$ has two
eigenvalues on the unit circle, under different restrictions, $\{
{\mathbf X}_t\}$ might be null recurrent (but not $\beta$-null
recurrent) or transient. If ${\mathbf A}$ has three or more eigenvalues
on the unit circle, the VAR(1) process $\{{\mathbf X}_t\}$ would be
transient, which indicates that the limit theory developed in this
paper would be not applicable.
\end{exa}

We next briefly discuss a nonlinear autoregressive model of the form:
%
\begin{equation}
\label{eq4.11} X_{t+1}=g(X_t,{\bolds\theta}_0)+e_{t+1},
\qquad t=1,2,\ldots,n.
\end{equation}
For this autoregression case, $\{e_t\}$ is not independent of $\{X_t\}$, and thus the proof strategy developed in this paper needs to be
modified. Following the argument in \citet{KT01}, in
order to develop an asymptotic theory for the parameter estimation in
the nonlinear autoregression (\ref{eq4.11}), we may need that the
process $\{X_t\}$ is Harris recurrent but not that the compound process
$\{(X_t,e_{t+1})\}$ is also Harris recurrent. This is because we
essentially have to consider sums of products like $\dot{g}(X_t,
{\bolds\theta}_0) e_{t+1} = \dot{g}(X_t, {\bolds\theta}_0)
(X_{t+1}-g(X_t,{\bolds\theta}_0))$, which are of the general form
treated in \citet{KT01}. The verification of the
Harris recurrence of $\{X_t\}$ has been discussed by \citet{L98} and
Example~\ref{exa21} given in Section~\ref{sec2.2} above. How to establish an asymptotic
theory for the parameter estimation of ${\bolds\theta}_0$ in model
(\ref{eq4.11}) will be studied in our future research.

\section{Simulated examples}\label{sec5}

In this section, we provide some simulation studies to compare the
finite sample performance of the proposed parametric estimation methods
and to illustrate the developed asymptotic theory.

\begin{exa}\label{exa51}
Consider the generalized linear model
defined by
%
\begin{equation}
\label{eq5.1} Y_t=\exp\bigl\{- \theta_0
X_t^2\bigr\}+e_t, \qquad
\theta_0=1, t=1,2,\ldots,n,
\end{equation}
where $\{X_t\}$ is generated by one of the three Markov processes:
\begin{longlist}[(iii)]
\item[(i)]  AR(1) process: $X_t=0.5X_{t-1}+x_t$,
\item[(ii)]  Random walk process: $X_t=X_{t-1}+x_t$,
\item[(iii)]  TAR(1) process: $X_t=0.5X_{t-1}I(|X_{t-1}|\leq
1)+X_{t-1}I(|X_{t-1}|>1)+x_t$,
$X_0=0$ and $\{x_t\}$ is a sequence of i.i.d. standard normal random
variables for the above three processes. The error process $\{e_t\}$ is
a sequence of i.i.d. ${\mathsf{N}}(0,0.5^2)$ random variables and
independent of $\{x_t\}$. In this simulation study, we compare the
finite sample behavior of the NLS estimator $\widehat{\theta}_n$ with
that of the MNLS estimator~$\overline{\theta}_n$, and the sample size
$n$ is chosen to be $500$, $1000$ and $2000$. The aim of this example
is to illustrate the asymptotic theory developed in Section~\ref
{sec3.1} as the regression function in (\ref{eq5.1}) is integrable when
$\theta_0>0$. Following the discussion in Section~\ref{sec2.2}, the
AR(1) process defined in (i) is positive recurrent, and the random
process defined in (ii) and the TAR(1) process defined in (iii) are
$1/2$-null recurrent.

We generate $500$ replicated samples for this simulation study, and
calculate the means and standard errors for both of the parametric
estimators in $500$ simulations. In the MNLS estimation procedure, we
choose $M_n=C_{\alpha}n^{1-\beta}$ with $\alpha=0.01$, where
$C_{\alpha
}$ is defined in (\ref{eq3.16}), $\beta=1$ for case (i), and $\beta
=1/2$ for cases (ii) and (iii). It is easy to find that $C_{0.01}=2.58$.
\end{longlist}

\begin{table}[t]
\caption{Means and standard errors for the estimators in Example~\protect\ref{exa51}}\label{tab51}
\begin{tabular*}{\textwidth}{@{\extracolsep{\fill}}lccc@{}}
\hline
\textbf{Sample size} & \textbf{500} & \textbf{1000} & \textbf{2000} \\
\hline
\multicolumn{4}{@{}c@{}}{The regressor $X_t$ is generated in case (i)}\\
NLS & 1.0036 (0.0481) & 1.0002 (0.0339) & 1.0000 (0.0245) \\
MNLS & 1.0036 (0.0481) & 1.0002 (0.0339) & 1.0000 (0.0245) \\[3pt]
\multicolumn{4}{@{}c@{}}{The regressor $X_t$ is generated in case (ii)}\\
NLS & 0.9881 (0.1783) & 0.9987 (0.1495) & 0.9926 (0.1393) \\
MNLS & 0.9881 (0.1783) & 0.9987 (0.1495) & 0.9926 (0.1393) \\[3pt]
\multicolumn{4}{@{}c@{}}{The regressor $X_t$ is generated in case (iii)}\\
NLS & 0.9975 (0.1692) & 1.0028 (0.1463) & 0.9940 (0.1301) \\
MNLS & 0.9975 (0.1692) & 1.0028 (0.1463) & 0.9940 (0.1301) \\
\hline
\end{tabular*}
\end{table}

The simulation results are reported in Table~\ref{tab51}, where the numbers in
the parentheses are the standard errors of the NLS (or MNLS) estimator
in the $500$ replications. From Table~\ref{tab51}, we have the following
interesting findings. (a) The parametric estimators perform better in
the stationary case (i) than in the nonstationary cases (ii) and (iii).
This is consistent with the asymptotic results obtained in Section~\ref
{sec3.1} such as Theorem~\ref{th3.1} and Corollaries \ref{cor31} and \ref{cor32}, which
indicate that the convergence rates of the parametric estimators can
achieve $O_P(n^{-1/2})$ in the stationary case, but only
$O_P(n^{-1/4})$ in the $1/2$-null recurrent case. (b) The finite sample
behavior of the MNLS estimator is the same as that of NLS estimator
since $\alpha=0.01$ means little sample information is lost. (c) Both
of the two parametric estimators improve as the sample size increases.
(d) In addition, for case (i), the ratio of the standard errors between
$500$ and $2000$ is $1.9633$ (close to the theoretical ratio $\sqrt
{4}=2$); for case (iii), the ratio of the standard errors between $500$
and $2000$ is $1.3005$ (close to the theoretical ratio
$4^{1/4}=1.4142$). Hence, this again confirms that our asymptotic
theory is valid.
\end{exa}

\begin{exa}\label{exa52}
Consider the quadratic regression model
defined by
%
\begin{equation}
\label{eq5.2} Y_t=\theta_0X_t^2+e_t,
\qquad \theta_0=0.5, t=1,2,\ldots,n,
\end{equation}
where $\{X_t\}$ is generated either by one of the three Markov
processes introduced in Example~\ref{exa51}, or by (iv) the unit root process:
\[
X_t=X_{t-1}+x_t, \qquad x_t=0.2x_{t-1}+v_t,
\]
in which $X_0=x_0=0$, $\{v_t\}$ is a sequence of i.i.d. ${\mathsf
{N}}(0,0.75)$ random variables, and the error process $\{e_t\}$ is
defined as in Example~\ref{exa51}. In this simulation study, we are interested
in the finite sample behavior of the MNLS estimator to illustrate the
asymptotic theory developed in Section~\ref{sec3.2} as the regression
function in (\ref{eq5.2}) is asymptotically homogeneous. For the
comparison purpose, we also investigate the finite sample behavior of
the NLS estimation, although we do not establish the related asymptotic
theory under the framework of null recurrent Markov chains. The sample
size $n$ is chosen to be $500$, $1000$ and $2000$ as in Example~\ref{exa51} and
the replication number is $R=500$. In the MNLS estimation procedure, as
in the previous example, we choose $M_n=2.58 n^{1-\beta}$, where
$\beta
=1$ for case (i), and $\beta=1/2$ for cases (ii)--(iv).
\end{exa}
\begin{table}[t]
\caption{Means and standard errors for the estimators in
Example~\protect\ref{exa52}}\label{tab52}
\begin{tabular*}{\textwidth}{@{\extracolsep{\fill}}lccc@{}}
\hline
\textbf{Sample size} & \textbf{500} & \textbf{1000} & \textbf{2000} \\
\hline
\multicolumn{4}{@{}c@{}}{The regressor $X_t$ is generated in case (i)}\\
NLS & 0.5002 (0.0095) & 0.4997 (0.0068) & 0.4998 (0.0050) \\
MNLS & 0.5003 (0.0126) & 0.4998 (0.0092) & 0.4997 (0.0064) \\[3pt]
\multicolumn{4}{@{}c@{}}{The regressor $X_t$ is generated in case (ii)}\\
NLS & 0.5000 ($2.4523\times10^{-4}$) & 0.5000 ($6.7110\times10^{-5}$) &
0.5000 ($2.7250\times10^{-5}$) \\
MNLS & 0.5000 ($2.4523\times10^{-4}$) & 0.5000 ($6.7112\times10^{-5}$)
& 0.5000 ($2.7251\times10^{-5}$) \\[3pt]
\multicolumn{4}{@{}c@{}}{The regressor $X_t$ is generated in case (iii)}\\
NLS & 0.5000 ($2.6095\times10^{-4}$) & 0.5000 ($8.4571\times10^{-5}$) &
0.5000 ($3.1268\times10^{-5}$) \\
MNLS & 0.5000 ($2.6095\times10^{-4}$) & 0.5000 ($8.4572\times10^{-5}$)
& 0.5000 ($3.1268\times10^{-5}$) \\[3pt]
\multicolumn{4}{@{}c@{}}{The regressor $X_t$ is generated in case (iv)}\\
NLS & 0.5000 ($2.1698\times10^{-4}$) & 0.5000 ($7.1500\times10^{-5}$) &
0.5000 ($2.6017\times10^{-5}$) \\
MNLS & 0.5000 ($2.1699\times10^{-4}$) & 0.5000 ($7.1504\times10^{-5}$)
& 0.5000 ($2.6017\times10^{-5}$) \\
\hline
\end{tabular*}
\end{table}

The simulation results are reported in Table~\ref{tab52}, from which, we have
the following conclusions. (a) For the regression model with
asymptotically homogeneous regression function, the parametric
estimators perform better in the nonstationary cases (ii)--(iv) than in
the stationary case (i). This finding is consistent with the asymptotic
results obtained in Sections~\ref{sec3.2} and \ref{sec4.2}. (b) The
MNLS estimator performs as well as the NLS estimator (in particular for
the nonstationary cases). Both the NLS and MNLS estimations improve as
the sample size increases.

\section{Empirical application}\label{sec6}

In this section, we give an empirical application of the proposed
parametric model and estimation methodology.

\begin{exa}\label{exa61}
Consider the logarithm of the UK to US export and
import data (in \pounds). These data come from the website:
\surl{https://www.\\ uktradeinfo.com/}, spanning from January 1996 to August 2013
monthly and with the sample size $n=212$. Let $X_t$ be defined as $\log
(E_t) + \log(p_t^{\mathrm{UK}}) - \log(p_t^{\mathrm{US}})$, where $\{E_t\}$ is the
monthly average of the nominal exchange rate, and $\{p_t^i\}$ denotes
the consumption price index of country $i$. In this example, we let $\{
Y_t\}$ denote the logarithm of either the export or the import value.

The data $X_t$ and $Y_t$ are plotted in Figures~\ref{fig61} and \ref{fig62},
respectively. Meanwhile, the real data application considered by \citet{GTY13} suggests that $\{X_t\}$ may follow the threshold
autoregressive model proposed in that paper, which is shown to be a
$1/2$-null recurrent Markov process. Furthermore, an application of the
estimation method by (\ref{eq3.5*}) gives $\beta_0 = 0.5044$. This
further supports that $\{X_t\}$ roughly follows a $\beta$-null
recurrent Markov chain with $\beta=1/2$.

\begin{figure}[t]

\includegraphics{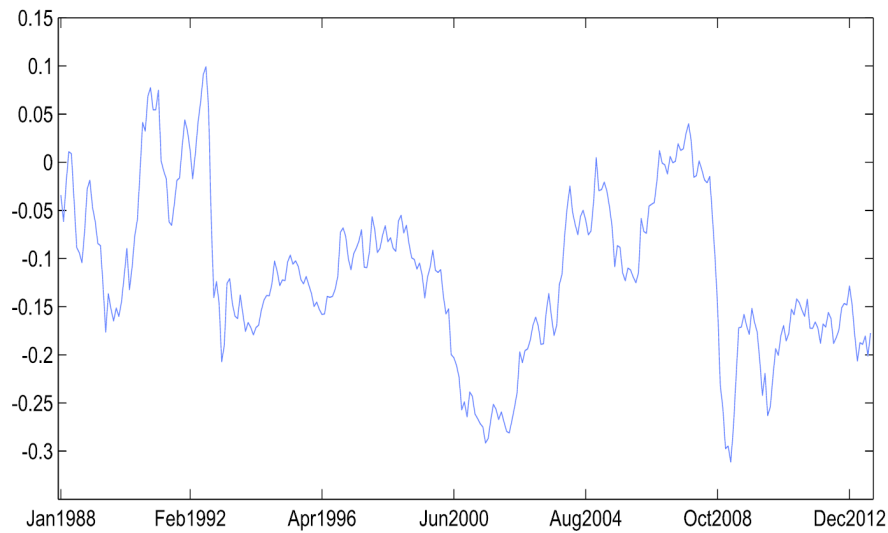}

\caption{Plot of the real exchange rate $X_t$.}\vspace*{-3pt}\label{fig61}
\end{figure}
\begin{figure}[b]

\includegraphics{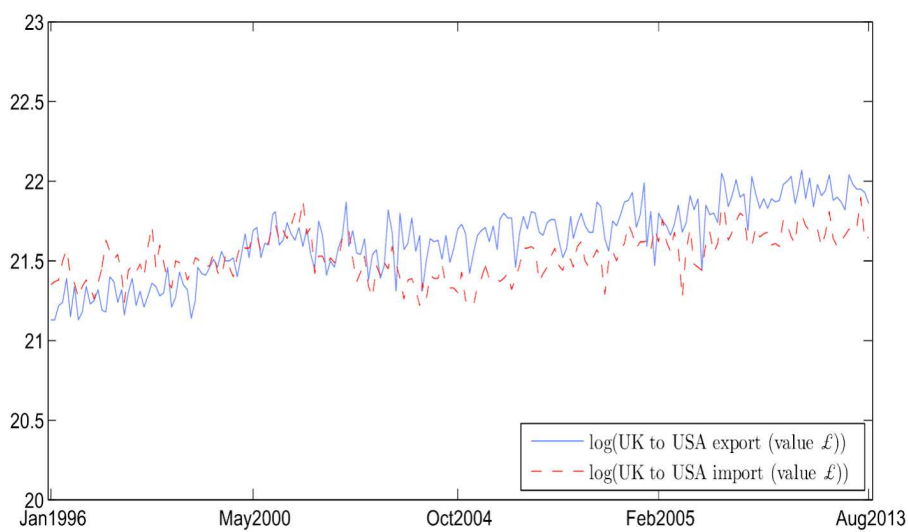}

\caption{Plot of the logarithm of the export and import
data $Y_t$.}\vspace*{-3pt}\label{fig62}
\end{figure}

To avoid possible confusion, let $Y_{\mathrm{ex},t}$ and $Y_{\mathrm{im},t}$ be the
export and import data, respectively. We are interested in estimating
the parametric relationship between $Y_{\mathrm{ex},t}$ (or $Y_{\mathrm{im},t}$) and
$X_t$. In order to find a suitable parametric relationship, we first
estimate the relationship nonparametrically based on $Y_{\mathrm{ex},t} =
m_{\mathrm{ex}}(X_t) + e_{t1}$ and $Y_{\mathrm{im},t}=m_{\mathrm{im}}(X_t)+e_{t2}$ [c.f., \citet{KMT07}], where $m_{\mathrm{ex}}(\cdot)$ and $m_{\mathrm{im}}(\cdot)$ are
estimated by
%
\begin{eqnarray}
\widehat{m}_{\mathrm{ex}}(x) &=& \frac{\sum_{t=1}^n K (({X_t -x})/{h}
) Y_{\mathrm{ex},t}}{\sum_{t=1}^n K (({X_t - x})/{h} )} \quad
\mbox{and}
\nonumber
\\[-8pt]
\label{eq6.1}
\\[-8pt]
\nonumber
 \widehat{m}_{\mathrm{im}}(x) &=& \frac{\sum_{t=1}^n K (({X_t -x})/{h} ) Y_{\mathrm{im},t}}{\sum_{t=1}^n K (({X_t - x})/{h} )},
\end{eqnarray}
where $K(\cdot)$ is the probability density function of the standard
normal distribution and the bandwidth $h$ is chosen by the conventional
leave-one-out cross-validation method. Then a parametric calibration
procedure (based on the preliminary nonparametric estimation) suggests
using a third-order polynomial relationship of the form
%
\begin{equation}
\label{eq6.2} Y_{\mathrm{ex}, t} = \theta_{\mathrm{ex},0} + \theta_{\mathrm{ex}, 1}
X_t + \theta_{\mathrm{ex}, 2} X_t^2 +
\theta_{\mathrm{ex}, 3} X_t^3+e_{\mathrm{ex},t}
\end{equation}
for the export data, where the estimated values (by using the method in
Section~\ref{sec3.2}) of $\theta_{\mathrm{ex}, 0}, \theta_{\mathrm{ex}, 1}, \theta_{\mathrm{ex},
2}$ and $\theta_{\mathrm{ex}, 3}$ are $21.666$, $5.9788$, $60.231$ and $139.36$,
respectively, and
%
\begin{equation}
\label{eq6.3} Y_{\mathrm{im}, t} = \theta_{\mathrm{im}, 0} + \theta_{\mathrm{im}, 1}
X_t + \theta_{\mathrm{im}, 2} X_t^2 +
\theta_{\mathrm{im}, 3} X_t^3+e_{\mathrm{im}, t}
\end{equation}
for the import data, where the estimated values of $\theta_{\mathrm{im},0},
\theta_{\mathrm{im},1}, \theta_{\mathrm{im},2}$ and $\theta_{\mathrm{im},3}$ are $21.614$,
$3.5304$, $37.789$ and $87.172$, respectively. Their plots are given in
Figures~\ref{fig63} and \ref{fig64}, respectively.

\begin{figure}[b]

\includegraphics{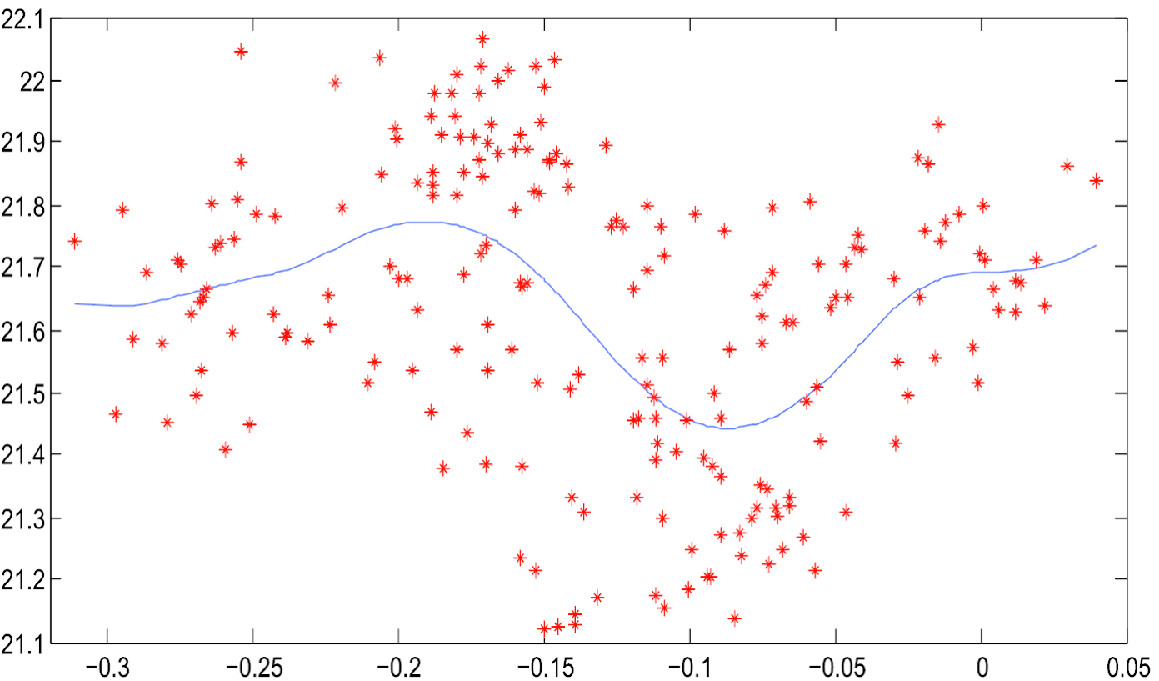}

\caption{Plot of the polynomial model fitting (\protect\ref{eq6.2}).}\label{fig63}
\end{figure}

\begin{figure}[t]

\includegraphics{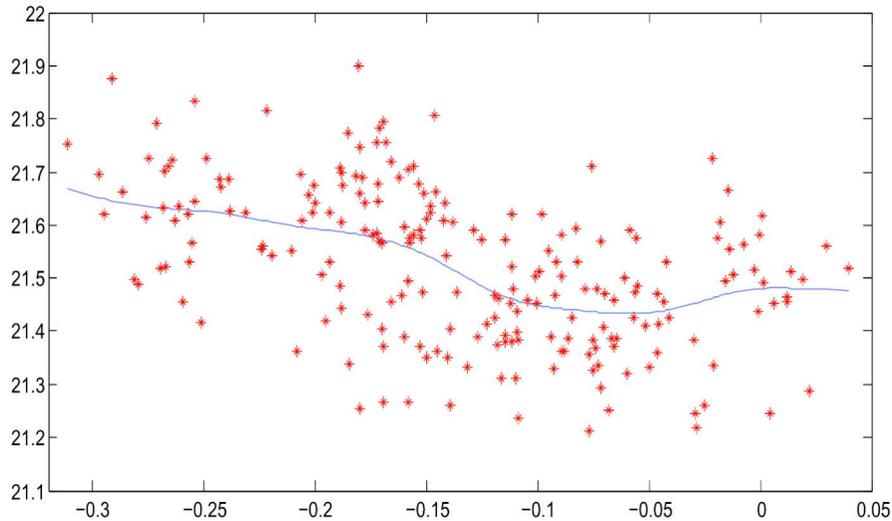}

\caption{Plot of the polynomial model fitting (\protect\ref{eq6.3}).}\label{fig64}
\end{figure}

While Figures~\ref{fig63} and \ref{fig64} suggest some relationship between the
exchange rate and either the export or the import variable, the true
relationship may also depend on some other macroeconomic variables,
such as, the real interest rate in the UK during the period. As
discussed in Section~4.3, we would like to extend the proposed models
from the univariate case to the multivariate case. As a future
application, we should be able to find a more accurate relationship
among the export or the import variable with the exchange rate and some
other macroeconomic variables.
\end{exa}

\section{Conclusions}\label{sec7}

In this paper, we have systematically studied the nonlinear regression
under the general Harris recurrent Markov chain framework, which
includes both the stationary and nonstationary cases. Note that the
nonstationary null recurrent process considered in this paper is under
Markov perspective, which, unlike PP, indicates that our methodology
has the potential of being extended to the nonlinear autoregressive
case. In this paper, we not only develop an asymptotic theory for the
NLS estimator of ${\bolds\theta}_0$ when $g(\cdot,\cdot)$ is
integrable, but also propose using a modified version of the
conventional NLS estimator for the asymptotically homogeneous $g(\cdot
,\cdot)$ and adopt a novel method to establish an asymptotic theory for
the proposed modified parametric estimator. Furthermore, by using the
log-transformation, we discuss the estimation of the parameter vector
in a conditional volatility function. We also apply our results to the
nonlinear regression with $I(1)$ processes which may be non-Markovian,
and establish an asymptotic distribution theory, which is comparable to
that obtained by PP. The simulation studies and empirical applications
have been provided to illustrate our approaches and results.

\begin{appendix}
\section*{Appendix A: Outline of the main proofs}\label{app}

In this Appendix, we outline the proofs of the main results in
Section~\ref{sec3}. The detailed proofs of these results are given in
Appendix~C of the supplemental document. The major difference between our proof
strategy and that based on the unit root framework [c.f., PP and
\citet{KR13}] is that our proofs rely on the limit
theorems for functions of the Harris recurrent Markov process (c.f.,
Lemmas \ref{lema1} and~\ref{lema2} below) whereas PP and \citet{KR13}'s
proofs use the limit theorems for integrated time series. We start with
two technical lemmas which are crucial for the  proofs of Theorems
\ref{th3.1} and \ref{th3.2}. The proofs for these two lemmas are given
in Appendix D of the supplemental document by \citet{LTG15}.

\renewcommand{\thelem}{A.1}
\begin{lem}\label{lema1}
Let $h_I(x,{\bolds\theta})$ be a
$d$-dimensional integrable function on $\Theta$ and suppose that
Assumption~\ref{ass31}\textup{(i)} is satisfied for $\{X_t\}$.
\begin{longlist}[(a)]
\item[(a)]  Uniformly for ${\bolds\theta}\in\Theta$, we have
%
\begin{equation}
\label{eqA.1} \frac{1}{N(n)}\sum_{t=1}^nh_I(X_t,{
\bolds\theta})=\int h_I(x,{\bolds\theta})\pi_s(dx)+o_P(1).
\end{equation}

\item[(b)]  If $\{e_t\}$ satisfies Assumption~\ref{ass31}\textup{(ii)}, we have, uniformly for
${\bolds\theta}\in\Theta$,
%
\begin{equation}
\label{eqA.2} \sum_{t=1}^nh_I(X_t,{
\bolds\theta})e_t=O_P \bigl(\sqrt {N(n)} \bigr).
\end{equation}
Furthermore, if $\int h_I(x,{\bolds\theta}_0)h_I^\tau
(x,{\bolds\theta}_0)\pi_s(dx)$ is positive definite, we have
%
\begin{equation}
\label{eqA.3} \qquad\frac{1}{\sqrt{N(n)}}\sum_{t=1}^nh_I(X_t,{
\bolds\theta }_0)e_t\stackrel{d}\longrightarrow{\sf N}
\biggl({\mathbf0}_d, \sigma ^2\int h_I(x,{
\bolds\theta}_0)h_I^\tau(x,{\bolds
\theta}_0)\pi _s(dx) \biggr).
\end{equation}
\end{longlist}
\end{lem}

\renewcommand{\thelem}{A.2}
\begin{lem}\label{lema2}
Let $h_{\mathrm{AH}}(x,{\bolds\theta})$ be
a $d$-dimensional asymptotically homogeneous function on $\Theta$ with
asymptotic order $\kappa(\cdot)$ (independent of ${\bolds\theta}$)
and limit homogeneous function $\overline{h}_{\mathrm{AH}}(\cdot,\cdot)$.
Suppose that $\{X_t\}$ is a null recurrent Markov process with the
invariant measure $\pi_s(\cdot)$ and Assumption~\ref{ass33}\textup{(iv)} are satisfied,
and $\overline{h}_{\mathrm{AH}}(\cdot,{\bolds\theta})$ is continuous on the
interval $[-1,1]$ for all ${\bolds\theta}\in\Theta$.
Furthermore, letting
\begin{eqnarray*}
{\bolds\Delta}_{\mathrm{AH}}(n,{\bolds\theta})&=&\sum
_{i=0}^{[M_n]}\overline{h}_{\mathrm{AH}}\biggl(
\frac{i}{M_n},{\bolds\theta }\biggr)\overline{h}_{\mathrm{AH}}^\tau
\biggl(\frac{i}{M_n},{\bolds\theta}\biggr)\pi _s
\bigl(B_{i+1}(1)\bigr)
\\[-2pt]
&&{}+\sum_{i=0}^{[M_n]}\overline{h}_{\mathrm{AH}}
\biggl(\frac{-i}{M_n},{\bolds \theta }\biggr)\overline{h}_{\mathrm{AH}}^\tau
\biggl(\frac{-i}{M_n},{\bolds\theta}\biggr)\pi _s
\bigl(B_{i+1}(2)\bigr)
\end{eqnarray*}
with $B_i(1)$ and $B_i(2)$ defined in Section~\ref{sec3.2}, there exist
a continuous and positive definite matrix ${\bolds\Delta
}_{\mathrm{AH}}({\bolds\theta})$ and a sequence of positive numbers $\{
\zeta
_{\mathrm{AH}}(M_n)\}$ such that $\zeta_0(M_n)/\zeta_{\mathrm{AH}}(M_n)$ is bounded,
$\zeta_0(l_n)/\zeta_0(M_n)=o(1)$ for $l_n\rightarrow\infty$ but
$l_n=o(M_n)$, and
\[
\lim_{n\rightarrow\infty}\frac{1}{\zeta_{\mathrm{AH}}(M_n)}{\bolds\Delta }_{\mathrm{AH}}(n,{
\bolds\theta})={\bolds\Delta}_{\mathrm{AH}}({\bolds \theta}),
\]
where $\zeta_0(\cdot)$ is defined in Section~\ref{sec3.2}.
\begin{longlist}[(a)]
\item[(a)]  Uniformly for ${\bolds\theta}\in\Theta$, we have
%
\begin{eqnarray}
&& \bigl[N(n){\mathbf J}_{\mathrm{AH}}(n,{\bolds\theta})
\bigr]^{-1}\sum_{t=1}^nh_{\mathrm{AH}}(X_t,{
\bolds\theta})h_{\mathrm{AH}}^\tau(X_t,{\bolds \theta})I
\bigl(|X_t|\leq M_n \bigr)
\nonumber
\\[-8pt]
\label{eqA.4}
\\[-8pt]
\nonumber
&&\qquad=I_d+o_P(1),
\end{eqnarray}
where ${\mathbf J}_{\mathrm{AH}}(n,{\bolds\theta})=\kappa^2(M_n)\zeta
_{\mathrm{AH}}(M_n){\bolds\Delta}_{\mathrm{AH}}({\bolds\theta})$.

\item[(b)]  If $\{e_t\}$ satisfies Assumption~\ref{ass31}\textup{(ii)}, we have, uniformly for
${\bolds\theta}\in\Theta$,
%
\begin{equation}
\label{eqA.5} {\mathbf J}_{\mathrm{AH}}^{-1/2}(n,{\bolds\theta})\sum
_{t=1}^nh_{\mathrm{AH}}(X_t,{
\bolds\theta})I \bigl(|X_t|\leq M_n \bigr)e_t=O_P
\bigl(\sqrt{N(n)} \bigr),
\end{equation}
and furthermore,
%
\begin{equation}
\label{eqA.6}\qquad N^{-1/2}(n){\mathbf J}_{\mathrm{AH}}^{-1/2}(n,{
\bolds\theta}_0)\sum_{t=1}^nh_{\mathrm{AH}}(X_t,{
\bolds\theta}_0)I \bigl(|X_t|\leq M_n
\bigr)e_t\stackrel{d}\longrightarrow{\sf N} \bigl({\mathbf0}_d,
\sigma ^2I_d \bigr).
\end{equation}
\end{longlist}
\end{lem}

\begin{pf*}{Proof of Theorem \protect\ref{th3.1}}
For Theorem~\ref{th3.1}(a), we only need to verify the following sufficient condition
for the weak consistency [\citet{J69}]: for a~sequence of positive
numbers $\{\lambda_n\}$,
%
\begin{equation}
\label{eqA.7} \frac{1}{\lambda_n} \bigl[L_{n,g}({\bolds\theta
})-L_{n,g}({\bolds\theta}_0) \bigr]=L^\ast({
\bolds\theta },{\bolds\theta}_0)+o_P(1)
\end{equation}
uniformly for ${\bolds\theta}\in\Theta$, where $L^\ast(\cdot
,{\bolds\theta}_0)$ is continuous and achieves a unique minimum at
${\bolds\theta}_0$. This sufficient condition can be proved by
using (\ref{eqA.1}) and (\ref{eqA.2}) in Lemma \ref{lema1}, and (\ref{eq3.3})
in Theorem~\ref{th3.1}(a) is thus proved. Combining the so-called
Cram\'
{e}r--Wold device in \citet{B68} and (\ref{eqA.3}) in Lemma
\ref{lema1}(b), we can complete the proof of the asymptotically normal
distribution in (\ref{eq3.4}). Details can be found in Appendix C of
the supplementary material.
\end{pf*}

\begin{pf*}{Proof of Corollary \protect\ref{cor31}}
The asymptotic distribution
(\ref{eq3.6}) can be proved by using (\ref{eq3.5}) and
Theorem~\ref{th3.1}(b).
\end{pf*}

\begin{pf*}{Proof of Corollary~\protect\ref{cor32}}
The convergence result (\ref
{eq3.7}) can be proved by using (\ref{eq2.3}) and (\ref{eq3.4}), and
following the proof of Lemma \ref{lema2} in \citet{GKLT14}. A detailed
proof is given in Appendix C of the supplementary material.
\end{pf*}

\begin{pf*}{Proof of Theorem  \protect\ref{th3.2}}
The proof is similar to
the proof of Theorem~\ref{th3.1} above. To prove the weak consistency,
similar to (\ref{eqA.7}), we need to verify the sufficient condition:
for a sequence of positive numbers $\{\lambda_n^\ast\}$,
%
\begin{equation}
\label{eqA.8} \frac{1}{\lambda_n^\ast} \bigl[Q_{n,g}({\bolds\theta
})-Q_{n,g}({\bolds\theta}_0) \bigr]=Q^\ast({
\bolds\theta },{\bolds\theta}_0)+o_P(1)
\end{equation}
uniformly for ${\bolds\theta}\in\Theta$, where $Q^\ast(\cdot
,{\bolds\theta}_0)$ is continuous and achieves a unique minimum at
${\bolds\theta}_0$. Using Assumption~\ref{ass33}(ii) and following the
proofs of (\ref{eqA.4}) and (\ref{eqA.5}) in Lemma \ref{lema2} (see Appendix D
in the supplementary material), we may prove (\ref{eqA.8}) and thus the
weak consistency result (\ref{eq3.10}). Combining the Cram\'{e}r--Wold
device and (\ref{eqA.6}) in Lemma \ref{lema2}(b), we can complete the proof of
the asymptotically normal distribution for $\overline{\bolds
\theta
}_n$ in (\ref{eq3.11}). More details are given in Appendix C of the
supplementary material.
\end{pf*}

\begin{pf*}{Proof of Corollary~\protect\ref{cor33}}
By using Theorem~\ref{th3.2}(b)
and (\ref{eq2.3}), and following the proof of Lemma \ref{lema2} in \citet{GKLT14},
we can directly prove (\ref{eq3.13}).
\end{pf*}

\begin{pf*}{Proof of Corollary~\protect\ref{cor34}}
The convergence result (\ref
{eq3.14}) follows from (\ref{eq3.13}) in Corollary~\ref{cor33} and (\ref
{eq3.15}) can be proved by using (\ref{eq3.11}) in Theorem~\ref
{th3.2}(b).
\end{pf*}
\end{appendix}

\section*{Acknowledgments}\label{sec8}
The authors are grateful to the Co-Editor, Professor Runze Li, an
Associate Editor and two referees for their valuable and constructive
comments and suggestions that substantially improved an earlier version
of the paper. Thanks also go to Professor Peter Phillips and other
colleagues who commented on this paper and the participants of various
conferences and seminars where earlier versions of this paper were
presented. This work was started when the first and third authors
visited the second author at Department of Mathematics, University of
Bergen in 2011.

\begin{supplement}[id=suppA]
\stitle{Supplement to ``Estimation in nonlinear regression with Harris
recurrent Markov chains''}
\slink[doi]{10.1214/15-AOS1379SUPP} 
\sdatatype{.pdf}
\sfilename{aos1379\_supp.pdf}
\sdescription{We provide some additional simulation
studies, the detailed proofs of the main results in Section~\ref{sec3},
the proofs of Lemmas \ref{lema1}
and \ref{lema2} and Theorem~\ref{th4.1}.}
\end{supplement}

%




\printaddresses
\end{document}